\documentclass[11pt,a4paper]{article}

\usepackage{amssymb} 
\usepackage{amsmath}
\usepackage{amsthm}
\usepackage[normalem]{ulem}
\usepackage{xcolor,cancel}

\usepackage{extarrows}

\usepackage[T1]{fontenc}
\usepackage[latin1]{inputenc}
\usepackage{color}

\usepackage{amsfonts}
\usepackage{amssymb}

\usepackage{empheq}

\usepackage[T1]{fontenc}   
\usepackage{lmodern}       
\usepackage{microtype}     
\usepackage{hyperref}      	
\usepackage{verbatim}	   	

\usepackage{geometry}
\usepackage{graphicx}
\usepackage{subfigure}
\usepackage{tikz}
\usepackage{xcolor}
\usepackage{placeins}
\usepackage{wrapfig}
\usepackage{caption}
\usepackage{setspace}
\usepackage{float}
\usepackage{ulem}

\hypersetup{
    colorlinks=true,
    linkcolor=blue,
    citecolor=blue,      
    urlcolor=blue
}

\swapnumbers
\numberwithin{equation}{section}

\theoremstyle{plain}
\newtheorem{theorem}{Theorem}[section]

\newtheorem{corollary}[theorem]{Corollary}

\theoremstyle{definition}

\newtheorem*{remark*}{Remark}

\DeclareMathOperator{\area}{area}
\DeclareMathOperator{\vol}{vol}
\DeclareMathOperator{\iso}{iso}

\DeclareMathOperator{\image}{image}
\DeclareMathOperator{\arsinh}{arsinh}
\DeclareMathOperator{\artanh}{artanh}

\DeclareMathOperator{\sn}{sn}

\DeclareMathOperator{\cd}{cd}
	
\title{Embedded Delaunay tori and their Willmore energy}

\author{Christian Scharrer
\thanks{Max Planck Institute for Mathematics, Vivatsgasse 7, 53111 Bonn, Germany.  Email: \texttt{Scharrer@mpim-bonn.mpg.de}.}}

\begin{document} 
\maketitle

\begin{abstract}
	A family of embedded rotationally symmetric tori in the Euclidean 3-space consisting of two opposite signed constant mean curvature surfaces that converge as varifolds to a double round sphere is constructed. Using complete elliptic integrals, it is shown that their Willmore energy lies strictly below $8\pi$. Combining such a strict inequality with previous works by \textsc{Keller--Mondino--Rivi\`ere} and \textsc{Mondino--Scharrer} allows to conclude that for every isoperimetric ratio there exists a smoothly embedded torus minimising the Willmore functional under isoperimetric constraint, thus completing the solution of the isoperimetric-constrained Willmore problem for tori. Similarly, we deduce the existence of smoothly embedded tori minimising the Helfrich functional with small spontaneous curvature. Moreover, it is shown that the tori degenerate in the moduli space which gives an application also to the conformally-constrained Willmore problem. Finally, because of their symmetry, the Delaunay tori can be used to construct spheres of high isoperimetric ratio, leading to an alternative proof of the known result for the genus zero case.
\end{abstract}

\section{Introduction} \label{sec:intro}
Given an immersed surface $f:\Sigma \to \mathbb R^3$, the Willmore functional $\mathcal W$ at $f$ is defined by
\begin{equation*}
	\mathcal W(f) = \int_\Sigma H^2\,\mathrm d\mu,
\end{equation*}
where the mean curvature $H$ is given by the arithmetic mean of the two principal curvatures, and $\mu$ is the Radon measure on $\Sigma$ corresponding to the pull back metric of the Euclidean metric in $\mathbb R^3$ along $f$. The isoperimetric ratio is defined by
\begin{equation}\label{eq:intro:def_iso}
	\iso(f) = 
	\begin{cases}
		\frac{\area(f)}{\vol(f)^\frac{2}{3}} & \text{if $\vol(f)>0$},\\
		\infty & \text{if $\vol(f) \leq 0$},
	\end{cases}
\end{equation}
where 
\begin{equation}\label{eq:algebraic_volume} 
	\area(f) = \int_{\Sigma}1\,\mathrm d\mu, \qquad \vol(f) = \frac{1}{3}\int_{\Sigma} n \cdot f \,\mathrm d \mu
\end{equation}
are the area and algebraic volume, and $n: \Sigma \to \mathbb S^2$ is the Gau{\ss} map.

The aim of this paper is to construct a family of embedded $C^{1,1}$-regular tori~$\mathbb T_{\mathrm D,c}$ in~$\mathbb R^3$ corresponding to $1<c<\gamma_0$ for some constant $\gamma_0>1$. The tori are rotationally symmetric and converge for $c\to1$ as varifolds to a round sphere of multiplicity 2 and radius 2. They are constructed out of two kinds of constant mean curvature surfaces that are well known in literature as \emph{Delaunay surfaces} (see~\cite{Delaunay1841}): The inner part of the tori has constant, strictly positive mean curvature; the outer part has constant, strictly negative mean curvature. These two pieces of Delaunay surfaces have matching normal vectors along the curve of intersection, leading to $C^{1,1}$-regularity of the patched surface. For a picture of the profile curve, see Figure~\ref{fig:d-torus}. The tori will be called \emph{Delaunay tori}. Their main property is stated in the following theorem which will be proven in Section~\ref{sec:d-tori}, using complete elliptic integrals. 
\begin{theorem} \label{thm:8pi-bound}
	There exists $\gamma_0>1$ such that the family of embedded Delaunay tori $\mathbb T_{\mathrm{D},c}$ corresponding to $1<c<\gamma_0$ satisfies
	\begin{equation*}
		\mathcal W(\mathbb T_{\mathrm D,c}) < 8\pi \qquad \text{whenever } 1< c< \gamma_0
	\end{equation*}
	and 
	\begin{equation*}
		\lim_{c \to 1+} \iso(\mathbb T_{\mathrm{D},c}) = \infty.
	\end{equation*}
\end{theorem}

In Section~\ref{sec:d-spheres}, it will be shown that the family of Delaunay tori $\mathbb T_{\mathrm D,c}$ can be used to construct a family of embedded spheres $\mathbb S_{\mathrm D,c}$ in $\mathbb R^3$ (for a picture of the profile curve, see Figure~\ref{fig:d-sphere}) with the following property. 
\begin{theorem}\label{thm:d-spheres}
	Let $\gamma_0>1$ be such that the Delaunay tori $\mathbb T_{\mathrm D,c}$ exist for all $1<c<\gamma_0$. Then, the family of spheres $\mathbb S_{\mathrm D,c}$ satisfies
	\begin{equation*}
	\mathcal W(\mathbb S_{\mathrm D,c}) = 4\pi + \frac{\mathcal W(\mathbb T_{\mathrm D,c})}{2} \qquad \text{for all $1<c<\gamma_0$}
	\end{equation*}
	as well as 
	\begin{equation*}
	\lim_{c\to1+}\iso(\mathbb S_{\mathrm D,c}) = \infty.
	\end{equation*}
\end{theorem}   
There are different notions of isoperimetric ratio in literature all of which are scaling invariant. The definition in this paper~\eqref{eq:intro:def_iso} is the same as in \cite{mondino2020strict} but differs from \cite{Schygulla} and \cite{MR3176354}. Nevertheless, it is easy to see that Theorem~\ref{thm:d-spheres} together with Theorem~\ref{thm:8pi-bound}
provides an alternative proof of the known result for spheres \cite[Lemma~1]{Schygulla}. 

Denote with $\mathcal S_1$ the space of smoothly immersed tori in $\mathbb R^3$. As a consequence of the Euclidean isoperimetric inequality, the isoperimetric ratio is minimised exactly by any parametrisation of a round sphere. Indeed, the round sphere is the only closed stable constant mean curvature surface~\cite{MR731682}. On the other hand, each smoothly embedded closed surface in $\mathbb R^3$ can be smoothly transformed arbitrarily close to a round sphere. This can be done using a one parameter family of M\"obius transformations whose centres of inversion approach a point on the surface (see \cite[Proposition~1]{scharrer:PhD} and \cite[Theorem~3.1]{yu2020uniqueness}). It follows that
\begin{equation*}
	\iso[\mathcal S_1] = (\sqrt[3]{36\pi},\infty]
\end{equation*}
where $\iso(\mathbb S^2) = \sqrt[3]{36\pi}$. The following main application will be proven in Section~\ref{sec:d-tori}; the proof will follow by combining Theorem~\ref{thm:8pi-bound} with previous works of  \textsc{Keller--Mondino--Rivi\`ere}~\cite{MR3176354} and \textsc{Mondino--Scharrer}~\cite{mondino2020strict}. 
\begin{corollary} \label{cor:existence_iso_constrained_tori}
	Let $\sqrt[3]{36\pi}<\sigma<\infty$. Then, there holds 
	\begin{equation}\label{eq:intro:8pi-bound}
		\beta_1(\sigma) := \inf\{\mathcal W(f) \mid f\in \mathcal S_1,\,\iso(f) = \sigma\}<8\pi
	\end{equation} 
	and the infimum in~\eqref{eq:intro:8pi-bound} is attained by a smoothly embedded minimiser $f_0 \in \mathcal S_1$.
\end{corollary} 
This completes the solution for the existence (and regularity) problem of isoperimetric constrained minimisers for the Willmore functional in the genus one case. The genus zero case was solved by \textsc{Schygulla}~\cite{Schygulla}. 
Notice that Corollary~\ref{cor:existence_iso_constrained_tori} is stated for all $\sigma \in \iso[\mathcal S_1]$ while Theorem~\ref{thm:8pi-bound} only holds for high isoperimetric ratios $\sigma$. In fact, the crucial part of the strict inequality~\ref{eq:intro:8pi-bound} is exactly that it holds true for high isoperimetric ratios. Indeed, the function $\beta_1$ is non-decreasing in $\sigma$, see \cite[Theorem~3.15]{scharrer:PhD} and \cite[Corollary~1.6]{mondino2020strict}.

Given an immersed surface $f:\Sigma \to \mathbb R^3$ and $c_0 \in \mathbb R$, the \emph{Helfrich functional} $\mathcal H_{c_0}$ at $f$ is defined by
\begin{equation*} 
\mathcal H_{c_0}(f) = \int_\Sigma (H - c_0)^2\,\mathrm d\mu.
\end{equation*}
In order to study the shape of lipid bilayer cell membranes, \textsc{Helfrich}~\cite{helfrich1973elastic} proposed the minimisation of $\mathcal H_{c_0}$ in the class of closed surfaces with given fixed area and given fixed volume. The constant $c_0$ is referred to as \emph{spontaneous curvature}. Existence of minimisers in the class of (possibly branched and bubbled) spheres was proven by \textsc{Mondino--Scharrer}~\cite{MR4076069}. Existence and regularity for minimisers with higher genus remains an open problem. Partial results were obtained by \textsc{Choksi--Veneroni}~\cite{MR3116014}, \textsc{Eichmann}~\cite{MR4129521}, and \textsc{Brazda--Lussardi--Stefanelli}~\cite{MR4098040}. Denote with $\mathcal S_0$ the space of smoothly immersed spheres in $\mathbb R^3$ and for all $\sigma > \sqrt[3]{36\pi}$ let
\begin{equation*}
	\beta_0(\sigma) = \inf\{\mathcal W(f) \mid f\in\mathcal S_0,\,\iso(f)=\sigma\}.
\end{equation*}
Suppose $A_0,V_0>0$ satisfy the isoperimetric inequality: $A_0^3 > 36\pi V_0^2$. Then, by Corollary~\ref{cor:existence_iso_constrained_tori} and \cite[Corollary~1.5]{mondino2020strict}, the following constant is strictly positive:
\begin{equation} \label{eq:Helfrich-epsilon}
	\varepsilon(A_0,V_0):= \frac{\sqrt{\min\{8\pi, 2\pi^2 + \beta_0(A_0/V_0^{2/3}) - 4\pi\}} - \sqrt{\beta_1(A_0/V_0^{2/3})}}{2\sqrt{A_0}} > 0.
\end{equation}
As another application of Theorem~\ref{thm:8pi-bound}, the following result on the existence of Helfrich tori will be proven in Section~\ref{sec:Helfrich_tori}.
\begin{corollary} \label{cor:Helfrich-tori}
	Suppose $A_0,V_0>0$ satisfy the isoperimetric inequality: $A_0^3 > 36\pi V_0^2$ and let $\varepsilon=\varepsilon(A_0,V_0)$ be defined as in Equation~\eqref{eq:Helfrich-epsilon}. Then, for each $c_0\in(-\varepsilon,\varepsilon)$, there exists a smoothly embedded torus $f_0 \in \mathcal S_1$ with 
	\begin{equation*}
		\area(f_0) = A_0, \qquad \vol(f_0) = V_0
	\end{equation*}
	and
	\begin{equation*}
		\mathcal H_{c_0}(f_0) = \inf\{\mathcal H_{c_0}(f) \mid f\in\mathcal S_1,\, \area(f)=A_0,\, \vol(f) = V_0\}.
	\end{equation*}
\end{corollary}

A further application of Theorem~\ref{thm:8pi-bound} can be found in the context of conformally constrained minimisation. We define the moduli space $\mathcal M$ for tori as a subset of the complex plane by
\begin{equation*}
	\mathcal M := \{x + iy \in\mathbb C \mid -1/2\leq x \leq 1/2,\,y>0,\,x^2 + y^2 \geq 1\}.
\end{equation*}
Any torus is conformally equivalent to a quotient $T^2_\omega:=\mathbb C/(\mathbb Z + \omega \mathbb Z)$ endowed with the Euclidean metric for some $\omega \in \mathcal M$, cf. \cite{MR1215481}.
Following \textsc{Ndiaye--Sch\"{a}tzle}~\cite{MR3403429}, we let
\begin{equation}\label{eq:intro:minimal_conf_energy}
	\mathcal M_{3,1}(\omega) := \inf\{\mathcal W(f) \mid f:T^2_\omega\to\mathbb R^3\text{ is conformal}\}.
\end{equation}
The following corollary can be found in \textsc{Ndiaye--Sch\"{a}tzle}~\cite[Proposition~D.1]{MR3403429}. In Section~\ref{sec:d-tori}, we will show how the Delaunay tori and Theorem~\ref{thm:8pi-bound} provide an alternative proof to the one in~\cite{MR3403429}.
\begin{corollary}\label{cor:conf-constrained}
	There exists a constant $y_0>1$ such that 
	\begin{equation}\label{eq:intro:conf_8pi-bound}
		\mathcal M_{3,1}(iy) < 8\pi \qquad \text{for $y\geq y_0$.}
	\end{equation}
\end{corollary} 
Recently, \textsc{Dall'Acqua--M{\"u}ller--Sch{\"a}tzle--Spener}~\cite{dall2020willmore} proved
that the Willmore flow of rotationally symmetric tori with initial energy at most $8\pi$ stays rotationally symmetric, exists for all times, and converges to the Clifford torus. The conformal class depends continuously on the time (see \cite[Proposition 4.2]{dall2020willmore}) while the Willmore energy is non-increasing in time.  
Consequently, since the conformal class of the Clifford torus is represented by the complex number $\omega = i$,
the function $\mathcal M_{3,1}(iy)$ is non-increasing in $y\geq1$. In particular, the strict inequality in~\eqref{eq:intro:conf_8pi-bound} becomes valid for all $y\geq1$. Hence, by \textsc{Kuwert--Sch\"atzle}~\cite{MR3024303}, the infimum $\mathcal M_{3,1}(iy)$ for $y\geq1$ in~\eqref{eq:intro:minimal_conf_energy} is attained by a smooth minimiser. For an existence result on the conformally constrained minimisation that does not rely on an $8\pi$-bound, see \textsc{Rivi\`ere}~\cite{MR3276154,MR3383803}. Explicit minimisers can be found in~\cite{MR3247390,MR3403429,MR4270047}.

In many classical problems related to the minimisation of the Willmore functional, strict $8\pi$-bounds such as in Equation~\eqref{eq:intro:8pi-bound} or Corollary~\ref{cor:conf-constrained} play a crucial role. One of the reasons is that by the \emph{Li--Yau inequality}~\cite{MR674407}, any immersed surface with Willmore energy strictly below $8\pi$ is actually embedded. Exemplary for the importance of $8\pi$-bounds is the Willmore flow. \textsc{Kuwert--Sch\"atzle}~\cite{MR2119722} showed that the Willmore flow of spheres exists for all time and converges to a round sphere provided the initial surface has Willmore energy less than $8\pi$. Later, \textsc{Blatt}~\cite{MR2591055} showed that this energy threshold is actually sharp. Recently, \textsc{Dall'Acqua--M{\"u}ller--Sch{\"a}tzle--Spener}~\cite{dall2020willmore} showed that the same energy threshold is also sharp for the long time existence of the Willmore flow of rotationally symmetric tori. Equally important, $8\pi$-bounds are needed in the direct method of the calculus of variations for the minimisation of the Willmore functional. This relates to the classical Willmore problem, the conformally constrained Willmore problem, and the isoperimetric constrained Willmore problem. In what follows, we illustrate the importance of $8\pi$-bounds for the Willmore functional in literature and point out potential directions for future research.

In the early 60s, \textsc{Willmore}~\cite{MR0202066} showed that the energy now bearing his name is bounded from below by $4\pi$ on the class of closed surfaces, with equality only for the round sphere. This inequality is sometimes referred to as \emph{Willmore's inequality}. More than two decades later, in an interesting work that connects Willmore surfaces with minimal surfaces, \textsc{Kusner}~\cite{MR996204} estimated the area of the celebrated \emph{Lawson surfaces}~\cite{MR0270280}. \textsc{Kusner's} work~\cite{MR996204} led to the $8\pi$-bound for the unconstrained minimal Willmore energy amongst surfaces of arbitrary genus. This was one of the key steps in proving existence and regularity of minimisers for the classical minimisation problem proposed by \textsc{Willmore} \cite{MR0202066}. Roughly speaking, the $8\pi$ bound prevents macroscopic bubbling in the direct method of calculus of variations, as by Willmore's inequality, each of the bubbles would cost at least $4\pi$ energy. Indeed, already \textsc{Simon}~\cite{MR1243525} used the $8\pi$-bound to obtain compactness (up to suitable M\"obius renormalisations) in the so called \emph{ambient approach}. Later, the $8\pi$-bound was also used in the \emph{parametric approach} to obtain compactness in the moduli space of higher genus surfaces by independent papers of \textsc{Kuwert--Li}~\cite{MR2928715} (building on top of previous work of \textsc{M\"uller--\v{S}ver\'{a}k}~\cite{MR1366547}) and \textsc{Rivi\`ere}~\cite{MR3008339} (building on top of \textsc{H\'elein's} moving frames technique~\cite{MR1913803}). 

Existence of smoothly embedded isoperimetric constrained Willmore spheres was proven by \textsc{Schygulla}~\cite{Schygulla}. Inspired by the computations of \textsc{Castro-Villarreal--Guven}~\cite{castro2007inverted}, \textsc{Schygulla}~\cite{Schygulla} applied a family of sphere inversions to a complete catenoid resulting in a family of closed surfaces with arbitrarily high isoperimetric ratios having one point of multiplicity two and Willmore energy exactly $8\pi$. Subsequently, they applied the Willmore flow for a short time around the point of multiplicity two, to obtain a family of surfaces with Willmore energy strictly below $8\pi$ and arbitrarily high isoperimetric ratios. This proved the strict $8\pi$-bound for isoperimetric constrained spheres. Moreover, \textsc{Schygulla}~\cite{Schygulla} showed that, as varifolds, isoperimetric constrained Willmore spheres (as well as the inverted catenoids) converge to a round sphere of multiplicity 2 as the isoperimetric ratio tends to infinity. Notice that the same holds true for the family of Delaunay tori constructed in this paper. His blow up result was analysed in more detail by \textsc{Kuwert--Li}~\cite{MR3842922}. They showed that any sequence of isoperimetric constrained Willmore spheres whose isoperimetric ratios diverge to infinity, indeed converges (up to subsequences, scaling, and translating) to two concentric round spheres of almost the same radii connected by a catenoidal neck. These kind of surfaces (i.e.\  two concentric round spheres of nearly the same radii connected by one or more catenoidal necks) were also constructed by means of cutting and paste techniques for different purposes in the works of \textsc{K\"uhnel--Pinkall}~\cite{MR868618}, \textsc{M\"uller--R\"oger}~\cite{MR3229052}, \textsc{Ndiaye--Sch\"{a}tzle}~\cite{MR3403429}, and \textsc{Wojtowytsch}~\cite{MR3744688}. One of the advantages of this technique is that it produces not only tori but surfaces of any topological type.

It is very natural to expect that, provided higher genus isoperimetric constrained minimisers exist, the blow up result of \textsc{Kuwert--Li}~\cite{MR3842922} can be generalised to the higher genus cases. To be more precise, we expect that any sequence of genus~$g$ isoperimetric constrained Willmore surfaces whose isoperimetric ratios diverge to infinity, converges (up to subsequences, scaling, and translating) to two concentric round spheres of nearly the same radii connected by $g + 1$ catenoidal necks. It is an interesting question whether or not the catenoidal necks in the limit have to be distributed over the double sphere in a certain way. In view of the solution presented in this paper, it is of course very tempting to conjecture that the catenoidal necks have to satisfy a balancing condition analogous to the one for constant mean curvature surfaces, see for instance \textsc{Kapouleas}~\cite{MR1100207,MR1317648}, or \textsc{Korevaar--Kusner--Solomon}~\cite{MR1010168}. That would mean that for tori, the two catenoidal necks necessarily end up being antipodal. \\ 

\textbf{Acknowledgements.}
The author was supported by the EPSRC as part of the MASDOC DTC at the University of Warwick, Grant No. EP/HO23364/1. Moreover, the author would like thank \textsc{Andrea Mondino}, \textsc{Filip Rindler}, and \textsc{Peter Topping} for hints and discussions on the subject. The author would also like to thank the referee for their careful reading of the original manuscript.

\section{Elliptic integrals} \label{sec:elliptic_integrals} 
\begin{figure}[H] 
	\centering
	\includegraphics[width=0.45\textwidth]{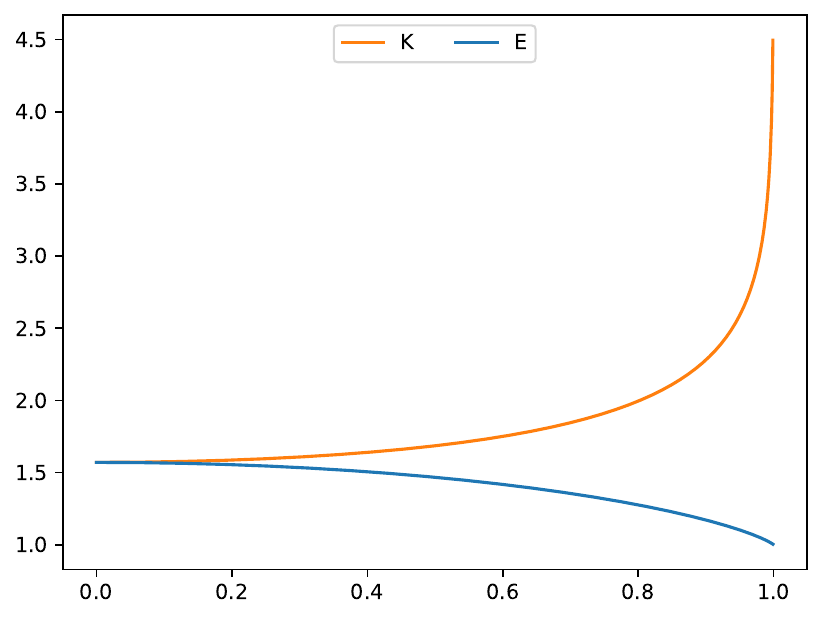}
	\caption{Complete elliptic integrals of the first and the second kind.}
\end{figure}
Elliptic integrals are functions defined as the value of common types of integrals that cannot be expressed in terms of elementary functions. They arise when computing geometric quantities such as the arc length of an ellipse or a hyperbola. In particular, they naturally occur in the context of constant mean curvature (Delaunay) surfaces of revolution. This is because the rotating curves of Delaunay surfaces are given by the roulette generated by ellipses and hyperbolas. In fact, all the quantities that are needed to construct the family of embedded Delaunay tori (see Section~\ref{sec:d-tori}) as well as their Willmore energy can be expressed in terms of \emph{complete elliptic integrals}. Given a so-called \emph{elliptic modulus} $k$, that is a real number $0<k<1$, the \emph{complete elliptic integral of the first kind} $K$ and the \emph{complete elliptic integral of the second kind} $E$ are defined by
\begin{equation*} 
	K(k) = \int_0^{\pi/2}\frac{\mathrm d\theta}{\sqrt{1 - k^2\sin^2(\theta)}}, \qquad E(k) = \int_0^{\pi/2} \sqrt{1 - k^2\sin^2(\theta)} \,\mathrm d \theta.
\end{equation*}   
All the formulas for elliptic integrals used in this paper can be found in the book of \textsc{Byrd--Friedman}~\cite{eInts}. The derivatives are given by
\begin{equation*} 
	\frac{\mathrm dK(k)}{\mathrm dk} = \frac{E(k)}{k(1-k^2)} - \frac{K(k)}{k}, \qquad \frac{\mathrm dE(k)}{\mathrm dk} = \frac{E(k) - K(k)}{k}.
\end{equation*}
The \emph{Gau{\ss} transformation} works as follows. Define the \emph{complementary modulus} $k'$ and the transformed modulus $k_1$ by
\begin{equation*}
	k' = \sqrt{1 - k^2}, \qquad k_1 = \frac{1 - k'}{1 + k'}.
\end{equation*}
Then, there holds (see~\cite[164.02]{eInts})
\begin{equation} \label{eq:eInts-Gauss_trans}
	K(k) = (1 + k_1)K(k_1), \qquad E(k) = (1 + k')E(k_1) - k'(1 + k_1)K(k_1).
\end{equation}
Moreover, $K$ grows like $\log(1/k')$, namely
\begin{equation} \label{eq:eInts-growth_of_K}
	\lim_{k \to 1-}\left( K(k) - \log(4/\sqrt{1 - k^2})\right) = 0
\end{equation}
and $E$ is bounded:
\begin{equation} \label{eq:eInts-boundedness_E}
	1 \leq E \leq \pi/2.
\end{equation}

\section{Surfaces of revolution}
A surface of revolution in $\mathbb R^3$ is given by a parametrisation $X$ of the type
\begin{equation*} \label{eq:surfaces_of_revolution} 
	X(t,\theta) = \left(f(t) \cos(\theta), f(t)\sin(\theta), g(t)\right)
\end{equation*} 
with parameters $t$ lying in an open interval and $0 \leq \theta \leq 2\pi$, where $f,g$ are real valued functions. The rotating curve $\gamma := (f,g)$ is referred to as \emph{meridian} or \emph{profile curve}. The underlying geometry is described by the coefficients of the first fundamental form
\begin{equation*}
	E^{I} = X_t\cdot X_t = {\dot f}^2 + {\dot g}^2 = |\dot \gamma|^2, \qquad F^{I}  = X_t \cdot X_\theta = 0, \qquad G^{I}  = X_\theta \cdot X_\theta = f^2
\end{equation*}
and the second fundamental form
\begin{equation*}
	L^{II}  = X_{tt} \cdot n = \frac{\dot f \ddot g - \ddot f \dot g}{|\dot \gamma|}, \qquad M^{II} = X_{t\theta} \cdot n = 0, \qquad N^{II} = X_{\theta\theta} \cdot n = \frac{f\dot g}{|\dot \gamma|} 
\end{equation*}
where the Gau{\ss} map $n$ is given by
\begin{equation*}
	n = \frac{X_t \times X_\theta}{|X_t \times X_\theta|}.
\end{equation*}
The mean curvature $H$ is defined as the arithmetic mean of the principal curvatures $\kappa_1,\kappa_2$, that is 
\begin{equation*}
	2H = \kappa_1 + \kappa_2 = \frac{L^{II}}{E^{I}} + \frac{N^{II}}{G^{I}} = \frac{\dot f \ddot g - \ddot f \dot g}{|\dot \gamma|^3} + \frac{\dot g}{f|\dot \gamma|}.
\end{equation*}
In this paper, we will focus on surfaces of revolution with constant mean curvature
\begin{equation} \label{eq:cmc_eqn}
	H = \frac{1}{2a}
\end{equation} 
for some given $0 \neq a \in \mathbb R$. These surfaces arise as critical points of the volume constrained area functional. Outside of a discrete set, one has $\dot g \neq 0$ and thus $\gamma(\varphi(t)) = (\rho(t),t)$ for some parameter transform $\varphi$ and some real valued function $\rho$. Hence, outside of a discrete set, Equation~\eqref{eq:cmc_eqn} can be turned into a second order ODE. Its solutions were first described by \textsc{Delaunay}~\cite{Delaunay1841} and are now named after him. More precisely, solutions for $a>0$ are called \emph{unduloids} and will be discussed in Section~\ref{sec:unduloids}; solutions for $a<0$ are called \emph{nodoids} and will be discussed in Section~\ref{sec:nodoids}. 

\section{Unduloids} \label{sec:unduloids}
\begin{figure}[H] 
	\centering
	\includegraphics[width=0.45\textwidth]{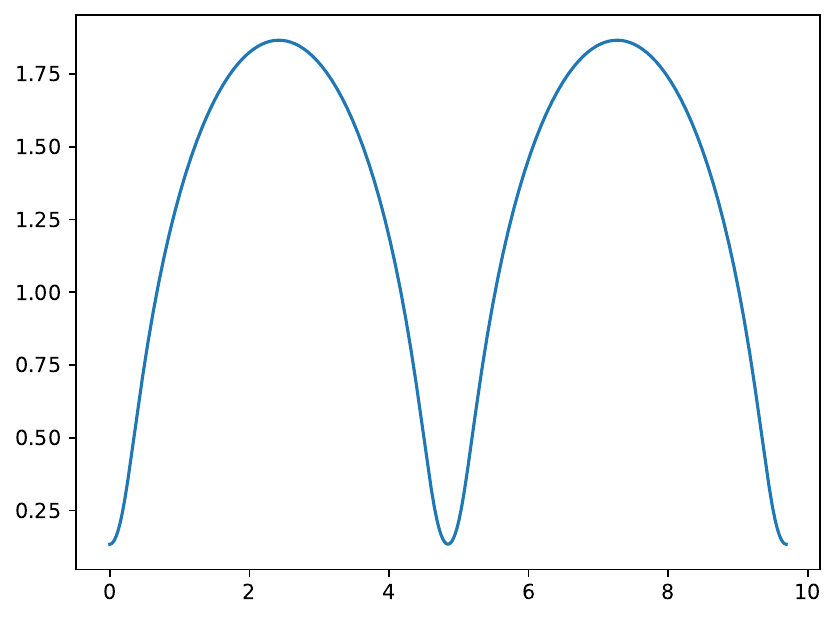}
	\caption{Profile curve of an unduloid with 
		$a=1$, $b=0.5$ 
		where the $f$-axis is vertical and the $g$-axis is horizontal.}
	\label{fig:unduloid-2periods}
\end{figure}
Unduloids are surfaces of revolution with constant, strictly positive mean curvature. Their rotating curve $(f,g)$ is given by the roulette of an ellipse with generating point being one of the foci. To be more precise, let $a>b>0$ and define $c = \sqrt{a^2 - b^2}$. Then, the equation
\begin{equation*}
	\frac{x^2}{a^2} + \frac{y^2}{b^2} = 1
\end{equation*}
describes a standard ellipse centred at the origin with width $2a$, height $2b$, and foci $\pm c$ on the $x$-axis. Rolling the ellipse without slipping along a line, each of the two focus points will describe a curve, \emph{the roulette}. 
This curve is \emph{periodic} in the following sense, where one period corresponds to one round of the ellipse: There exist a parametrisation $(f,g)$ of the roulette as well as $\lambda>0$ and $L\in\mathbb R$ referred to as \emph{period length} and \emph{extrinsic period length} such that
\begin{equation}\label{eq:periodic}
	\begin{split}
		f(t + k\lambda) & = f(t) \\
		g(t + k\lambda) & = g(t) + kL
 	\end{split}
\end{equation}  
for all $t\in\mathbb R$ and all $k\in\mathbb Z$. Indeed, the extrinsic length $L$ is given by the circumference of the ellipse, see~\eqref{eq:unduloid-extrinsic_length}. Figure~\ref{fig:unduloid-2periods} shows two periods of the roulette using the following parametrisation found by \textsc{Bendito--Bowick--Medina}~\cite{BBM2014} with period length $\lambda = 2\pi$:
\begin{align} \label{eq:unduloid-f_component}
	f(t) & = b\frac{a - c \cos(t)}{\sqrt{a^2 - c^2 \cos^2(t)}} \\ \label{eq:unduloid-g_component}
	g(t) & = \int_0^t\sqrt{a^2 - c^2 \cos^2(x)}\,\mathrm dx - c \sin(t)\frac{a - c \cos(t)}{\sqrt{a^2 - c^2 \cos^2(t)}}
\end{align}
for $0\leq t \leq 2\pi$ with coefficients of the first fundamental form
\begin{equation*} 
	E^{I} = \frac{a^2b^2}{(a + c\cos(t))^2}, \qquad G^{I} = b^2\frac{a - c \cos(t)}{a + c \cos(t)},
\end{equation*} 
speed
\begin{equation}\label{eq:speed-unduloid}
	|(f,g)'(t)| = \frac{ab}{a+c\cos(t)},
\end{equation}
and mean curvature
\begin{equation} \label{eq:unduloid-mean_curvature}
H = \frac{1}{2a}.
\end{equation}
The extrema are given by
\begin{equation} \label{eq:unduloid-extrem_values}
	\min \image f = f(0) = a - c, \qquad \max \image f = f(\pi) = a+c.
\end{equation}

\subsection{Extrinsic period length}
We compute the extrinsic length $L$ of one period. Using the parametrisation~\eqref{eq:unduloid-g_component}, we see
\begin{align*}
	L & = |g(2\pi) - g(0)| = a\int_0^{2\pi}\sqrt{1 - \frac{c^2}{a^2}\cos^2(x)} \,\mathrm dx = 4a\int_0^{\pi/2}\sqrt{\Bigl(1 - \frac{c^2}{a^2}\Bigr) + \frac{c^2}{a^2}\sin^2(x)} \,\mathrm dx \\
	& = 4a\sqrt{1 - \frac{c^2}{a^2}} \int_0^{\pi/2}\sqrt{1 + n^2\sin^2(x)}\,\mathrm dx
\end{align*}
for $n = c/b$. Letting $k^2 = n^2/(1 + n^2)$, (282.03) and (315.02) in \cite{eInts} imply 
\begin{equation*}
	\int_0^{\pi/2}\sqrt{1 + n^2\sin^2(x)}\,\mathrm dx = \frac{1}{k'}E(k),
\end{equation*}
where we used the special values $\sn(0)=0$ and $\cd(K)=0$ taken from (122.01), (122.02), and (120.02) in \cite{eInts}.
Thus, since
\begin{equation*}
	k = \frac{c}{a}, \qquad k'=\sqrt{1 - \frac{c^2}{a^2}}, 
\end{equation*}
it follows that
\begin{equation}\label{eq:unduloid-extrinsic_length}
	L = 4aE(k); \qquad k = \frac{c}{a}. 
\end{equation}

\subsection{Area computation}
Next, a formula for the area will be determined. Using the Weierstra{\ss} substitution, the area $\mathcal A$ of the rotational symmetric surface corresponding to one period can be computed by
\begin{align*}
	&\int_0^{2\pi}\int_0^{2\pi}\sqrt{EG}\,\mathrm d\theta\,\mathrm dt = 2\pi ab^2 \int_0^{2\pi} \sqrt{\frac{a - c\cos(t)}{(a+c\cos(t))^3}}\,\mathrm dt = 4\pi ab^2 \int_0^{\pi} \sqrt{\frac{a + c\cos(t)}{(a-c\cos(t))^3}}\,\mathrm dt \\
	&\quad = 4\pi ab^2\int_0^\infty\sqrt{\frac{a + c \frac{1-x^2}{1+x^2}}{(a - c \frac{1-x^2}{1+x^2})^3}} \frac{2\mathrm dx}{1 + x^2} = 8\pi ab^2\int_0^\infty\sqrt{\frac{a(1+x^2) + c(1-x^2)}{\bigl(a(1+x^2) - c(1 -x^2)\bigr)^3}}\,\mathrm dx \\
	&\quad = 8\pi ab^2\int_0^\infty \sqrt{\frac{(a + c) + (a - c)x^2}{\bigl((a - c) + (a + c)x^2\bigr)^3}}\,\mathrm dx = 8\pi ab^2\sqrt{\frac{a-c}{(a+c)^3}}\int_0^\infty\sqrt{\frac{\tilde a^2 + t^2}{(\tilde b^2 + t^2)^3}}\,\mathrm dt 
\end{align*}
for $\tilde a^2 = (a+c)/(a-c)$ and $\tilde b^2 = (a-c)/(a+c)$. The last integral can be transformed into a complete elliptic integral of the second kind using \cite[221.01]{eInts} with $k^2 = 1 - \tilde b^2/\tilde a^2$ and $g = 1/\tilde a$:
\begin{equation*}
	\int_0^\infty\sqrt{\frac{\tilde a^2 + t^2}{(\tilde b^2 + t^2)^3}}\,\mathrm dt = \frac{g}{k'^2}E(k).
\end{equation*}
One can show that 
\begin{equation*}
	k^2 = \frac{4ac}{(a+c)^2}, \qquad k' = \frac{a - c}{a+c}, \qquad \frac{g}{k'^2} = \frac{(a+c)^{3/2}}{(a - c)^{3/2}}
\end{equation*}
and thus 
\begin{equation} \label{eq:unduloid-area}
	\mathcal A = 8\pi a(a+c)E(k);\qquad k=\frac{2\sqrt{ac}}{a+c}.
\end{equation}
Notice that this coincides with the area formula for unduloids computed for a different parametrisation in~\cite{MR2381785} and~\cite{MR1977569}. 

\subsection{Volume computation}\label{sec:volume-unduloid}
One can show that for $g$ as defined in~\eqref{eq:unduloid-g_component} there holds
\begin{equation*}
	g'(t) = ab^2\frac{a - c\cos(t)}{(a^2 - c^2\cos^2(t))^{\frac{3}2}}.
\end{equation*}
Hence, using the Weierstra{\ss} substitution, the volume $\mathcal V$ of the rotational symmetric surface corresponding to one period can be computed as
\begin{align*}
	\mathcal V & = \pi\int_0^{2\pi}f^2g'\,\mathrm dt = 2\pi ab^4\int_0^\pi \sqrt{\frac{a + c\cos(t)}{(a - c\cos(t))^5}}\,\mathrm dt = 4\pi ab^4\int_0^\infty\sqrt{\frac{a + c\frac{1-x^2}{1+x^2}}{(a-c\frac{1-x^2}{1+x^2})^5}}\frac{\mathrm dx}{1+x^2}\\
	&=4\pi ab^4\int_0^\infty \sqrt{\frac{a+c + (a-c)x^2}{(a-c +(a+c)x^2)^5}}(1+x^2)\,\mathrm dx\\
	&=4\pi ab^4\int_0^{\frac{\pi}{2}}\sqrt{\frac{a+c + (a-c)\frac{a-c}{a+c}\tan^2(s)}{(a-c +(a+c)\frac{a-c}{a+c}\tan^2(s))^5}}\left(1+\frac{a-c}{a+c}\tan^2(s)\right)\sqrt{\frac{a-c}{a+c}}\frac{\mathrm ds}{\cos^2(s)}\\
	&=4\pi a\int_0^{\frac{\pi}{2}}\sqrt{(a+c)^2\cos^2(s)+(a-c)^2\sin^2(s)}\left((a+c)\cos^2(s) + (a-c)\sin^2(s)\right)\,\mathrm ds.
\end{align*}
The last line can be expressed in terms of complete elliptic integrals. For the purpose of this paper it is enough to notice that
\begin{equation}\label{eq:volume_convergence-unduloid}
	\mathcal V \xrightarrow{b\downarrow0}16\pi a^3\int_0^{\frac{\pi}{2}}\cos^3(s)\,\mathrm ds = \frac{4\pi}{3}(2a)^3.
\end{equation}
	
\section{Nodoids} \label{sec:nodoids}
\begin{figure}[H]
	\begin{minipage}[b]{0.45\linewidth}
		\centering
		\includegraphics[width=\textwidth]{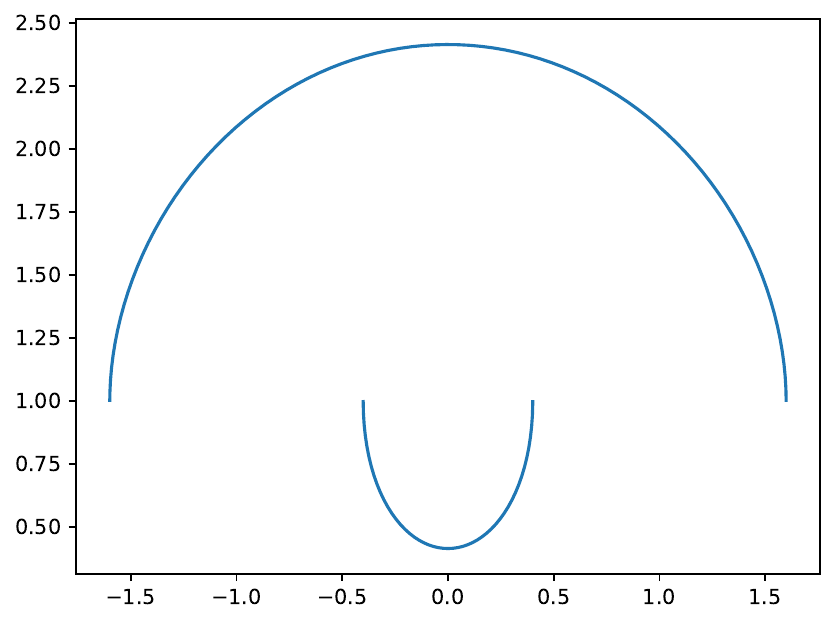}
		\caption{Separate roulettes $(f_\pm,g_\pm)$ (bottom/top) for $a=b=1$.}
		\label{fig:nodoid-2branches}
	\end{minipage}
	\hspace{0.5cm}
	\begin{minipage}[b]{0.45\linewidth}
		\centering
		\includegraphics[width=\textwidth]{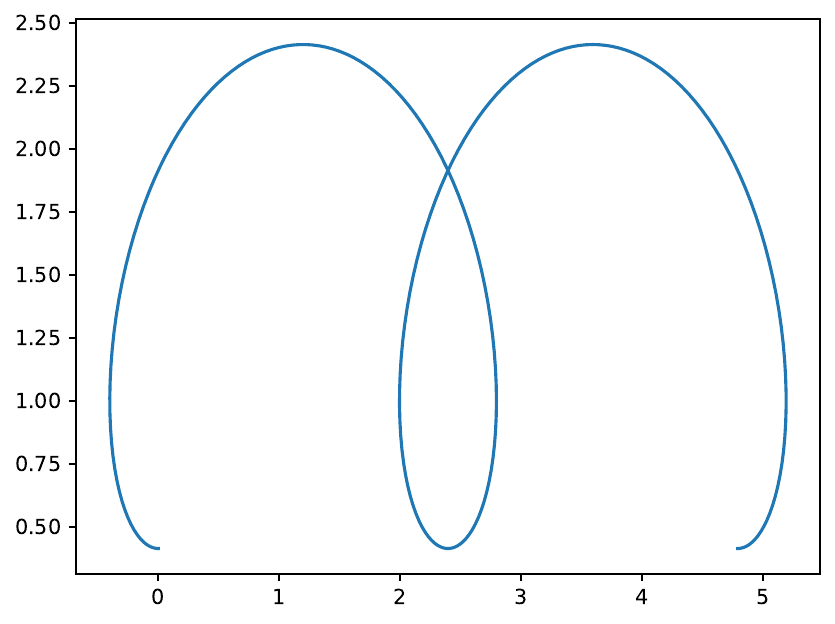}
		\caption{Both roulettes patched together for 2 periods and $a=b=1$.}
		\label{fig:nodoid-2periods}
	\end{minipage}
\end{figure}
Nodoids are surfaces of revolution with constant, strictly negative mean curvature. Their rotating curve $(f,g)$ is given by the roulette of a hyperbola with generating points given by the foci. To be more precise, let $a,b>0$ and define $c = \sqrt{a^2 + b^2}$. Then, the equation
\begin{equation*}
	\frac{x^2}{a^2} - \frac{y^2}{b^2} = 1
\end{equation*}
describes a hyperbola in canonical form with distance $a$ to the centre and foci $\pm c$ on the $x$-axis. Rolling the right branch of the hyperbola without slipping along a line, each of the two focus points will describe a curve, \emph{the roulette}. \textsc{Bendito--Bowick--Medina}~\cite{BBM2014} found parametrisations $(f_{\pm},g_{\pm})$ of the roulettes, where $(f_+,g_+)$ corresponds to the focus $(c,0)$ and $(f_-,g_-)$ is the reflected roulette corresponding to the focus $(-c,0)$:
\begin{align} \label{eq:nodoid-f_component}
	f_\pm(t) & = b\frac{c \cosh(t) \mp a}{\sqrt{c^2 \cosh^2(t) - a^2}} \\ 		\label{eq:nodoid-g_component}
	g_\pm(t) & = \int_0^t\sqrt{c^2 \cosh^2(x) - a^2}\,\mathrm dx - c \sinh(t)\frac{c \cosh(t)\mp a}{\sqrt{c^2 \cosh^2(t) - a^2}}
\end{align}
with parameter $t$ running through all of $\mathbb R$, coefficients of the first fundamental form
\begin{equation*} 
	E_\pm = \frac{a^2b^2}{(c\cosh(t) \pm a)^2}, \qquad G_\pm = b^2\frac{c \cosh(t)\mp a}{c \cosh(t)\pm a},
\end{equation*} 
speed
\begin{equation}\label{eq:nodoid-speed}
	|(f_\pm,g_\pm)'(t)| = \frac{ab}{c\cosh(t)\pm a},
\end{equation}
and mean curvature
\begin{equation} \label{eq:nodoid-mean_curvature}
	H_\pm = -\frac{1}{2a}.
\end{equation}
One has 
\begin{equation} \label{eq:nodoid-extrem_values}
	\begin{split}
		\max \image f_+ = \lim_{t\to \pm\infty} f_+(t) = b,\qquad \min \image f_+ = f_+(0) = c - a \\
		\min \image f_- = \lim_{t \to \pm \infty} f_-(t) = b \qquad \max \image f_- = f_-(0) = c+a
	\end{split} 
\end{equation}
and thus, after translation along the axis of rotation, the two roulettes corresponding to the foci $(\pm c,0)$ can be glued together into one 
curve. After suitable reparametrisation, this curve can be extended to a \emph{periodic} curve, i.e.\ a curve that satisfies~\eqref{eq:periodic}, see Figure~\ref{fig:nodoid-2periods}. Analysing the ODE~\eqref{eq:nodoid-mean_curvature}, one can show that the curve obtained by glueing is $C^2$-regular, see Equation~(3.3), Proposition~3, and Proposition~5 in \cite{MR1807955}. Indeed, in Section~\ref{subsec:ex_length_nod}, we will give a smooth parametrisation of the glued curve.

\subsection{Extrinsic period length} \label{subsec:ex_length_nod}
In order to determine the extrinsic length $L$ of one period, we perform the parameter transformation $t = \arsinh(u)$ and $u = \frac{b}{c}\tan(s)$ for $-\frac{\pi}{2} < s < \frac{\pi}{2}$. We have
\begin{align}\nonumber
	f_\pm(t) & = b\frac{c\sqrt{1 + u^2} \mp a}{\sqrt{c^2(1+u^2) - a^2}}=\frac{c\sqrt{1+u^2}\mp a}{\sqrt{1 + \frac{c^2}{b^2}u^2}} = \frac{c\sqrt{1 + \frac{b^2}{c^2}\tan^2(s)}\mp a}{\sqrt{1 + \tan^2(s)}}\\ \label{eq:nodoid-transformed_f}
	&=\sqrt{c^2\cos^2(s) +b^2\sin^2(s)} \mp a \cos(s)
\end{align}
and
\begin{align*}
	&\int_0^{t}\sqrt{c^2\cosh^2(x) - a^2}\,\mathrm dx = \int_0^u\sqrt{\frac{c^2(1+x^2) - a^2}{1 + x^2}}\,\mathrm dx = \frac{b}{c}\int_0^{\tan(s)}\sqrt{\frac{c^2\frac{b^2}{c^2}x^2 + b^2}{1+\frac{b^2}{c^2}x^2}}\,\mathrm dx \\
	&\quad=b^2\int_0^s\sqrt{\frac{1 + \tan^2(x)}{c^2 + b^2\tan^2(x)}}\tan'(x)\,\mathrm dx = b^2\int_0^s\frac{\tan'(x)\,\mathrm dx}{\sqrt{c^2\cos^2(x) +b^2\sin^2(x)}} \\
	&\quad=-b^2\int_0^s\frac{a^2\sin(x)\cos(x)}{(c^2-a^2\sin^2(x))^{\frac{3}{2}}}\tan(x)\,\mathrm dx + b^2\frac{\tan(s)}{\sqrt{b^2+a^2\cos^2(s)}}.
\end{align*}
Moreover, there holds $c\sinh(t) = b\tan(s)$ and 
\begin{align*}
	&\tan(s)\left[\frac{b^2}{\sqrt{b^2+a^2\cos^2(s)}} - \sqrt{b^2 + a^2\cos^2(s)} \pm a\cos(s)\right] \\
	&\quad = -a^2\frac{\sin(s)\cos(s)}{\sqrt{b^2+a^2\cos^2(s)}} \pm a\sin(s).
\end{align*}
It follows
\begin{equation} \label{eq:nodoid-transfomred_g}
	g_\pm(s) = -a^2b^2\int_0^s\frac{\sin^2(x)\,\mathrm dx}{(c^2 - a^2\sin^2(x))^\frac{3}{2}} -a^2\frac{\sin(s)\cos(s)}{\sqrt{b^2+a^2\cos^2(s)}} \pm a\sin(s).
\end{equation}
Define $k^2 = a^2/c^2$ and	
\begin{equation}\label{eq:ext_period_length_nodoid}
	L = \frac{4k^2b^2}{c}\int_0^{\pi/2}\frac{\sin^2(x)\,\mathrm dx}{(1 - k^2\sin^2(x))^{\frac{3}{2}}}.
\end{equation}
Extending the definition of the functions $f_\pm$ and $g_\pm$ in \eqref{eq:nodoid-transformed_f} and \eqref{eq:nodoid-transfomred_g} for all $s\in\mathbb R$, one readily verifies
\begin{equation}\label{eq:periodic-nodoid}
	f_{\pm}(s+\pi) = f_{\mp}(s), \qquad g_{\pm}(s+\pi) = -\frac{L}{2} + g_{\mp}(s)
\end{equation}
and thus
\begin{equation*}
	f_{\pm}(s+2\pi) = f_{\pm}(s), \qquad g_{\pm}(s+2\pi) = -L + g_{\pm}(s).
\end{equation*}
Hence, we see that both curves $\mathbb R \to \mathbb R^2$ with $s\mapsto(f_+(s),g_+(s))$ and $s\mapsto(f_-(s),g_-(s))$ parametrise the whole periodic curve resulting from the patched roulettes in Figure~\ref{fig:nodoid-2periods} with period length $\lambda = 2\pi$ and extrinsic period length $-L$ given in~\eqref{eq:ext_period_length_nodoid}, cf.\ \eqref{eq:periodic}. These curves are smooth. Using (281.05), (318.02) from \cite{eInts} as well as the special values $\sn(0) = \cd(K) = 0$, we find
\begin{equation} \label{eq:nodoid-extrinsic_length}
	L = 4c\left[E(k) - k'^2K(k)\right]; \qquad k = \frac{a}{c}.
\end{equation}

\subsection{Area computation}
Next, the area $\mathcal A$ of the rotational symmetric surface corresponding to one period will be computed. Using the parameter transformations $t = \artanh(x)$ and $x = \sin(t)$, one infers
\begin{align*}
	&\int_{-\infty}^{\infty}\sqrt{E_+G_+}\,\mathrm dt = ab^2\int_{-\infty}^{\infty}\sqrt{\frac{c\cosh(t) - a}{(c\cosh(t) + a)^3}}\,\mathrm dt = 2ab^2\int_0^1\sqrt{\frac{c\frac{1}{\sqrt{1-x^2}} - a}{(c\frac{1}{\sqrt{1-x^2}} + a)^3}} \frac{\mathrm dx}{1 - x^2}\\
	& \quad = 2ab^2\int_0^1 \sqrt{\frac{c - a\sqrt{1-x^2}}{(c + a\sqrt{1-x^2})^3}}\frac{\mathrm dx}{\sqrt{1-x^2}} = 2ab^2\int_0^{\pi/2}\sqrt{\frac{c - a\cos(t)}{(c + a \cos(t))^3}}\,\mathrm dt \\
	& \quad = 2ab^2\int_{\pi/2}^\pi\sqrt{\frac{c + a\cos(t)}{(c - a \cos(t))^3}}\,\mathrm dt
\end{align*}
and similarly,
\begin{equation*}
	\int_{-\infty}^{\infty}\sqrt{E_-G_-}\,\mathrm dt = 2ab^2\int_0^{\pi/2}\sqrt{\frac{c + a\cos(t)}{(c - a \cos(t))^3}}\,\mathrm dt.
\end{equation*}
Consequently,
\begin{equation*}
	\mathcal A = \int_0^{2\pi}\int_{-\infty}^{\infty}\sqrt{E_+G_+} + \sqrt{E_-G_-}\,\mathrm dt\,\mathrm d\theta = 4\pi ab^2\int_{0}^\pi\sqrt{\frac{c + a\cos(t)}{(c - a \cos(t))^3}}\,\mathrm dt.
\end{equation*}
In Section~\ref{sec:unduloids} it was shown that
\begin{equation*}
	\int_{0}^\pi\sqrt{\frac{c + a\cos(t)}{(c - a \cos(t))^3}}\,\mathrm dt = \frac{2E(k)}{c - a}; \qquad k = \frac{2\sqrt{ac}}{a+c}.
\end{equation*} 
Thus, it follows that
\begin{equation} \label{eq:nodoid-area}
	\mathcal A = 8\pi a(a+c)E(k);\qquad k = \frac{2\sqrt{ac}}{a+c}.
\end{equation}

\subsection{Volume computation}\label{sec:volume-nodoid}
Using the parametrisation in \eqref{eq:nodoid-transformed_f} and \eqref{eq:nodoid-transfomred_g}, we find
\begin{align*}
	f_\pm(s)^2& =2a^2\cos^2(s) + b^2\mp 2a\cos(s)\sqrt{a^2\cos^2(s) +b^2} ,\\
	g_\pm'(s) & = -\frac{a^2\cos^2(s)}{\sqrt{a^2\cos^2(s)+b^2}} \pm a\cos(s).
\end{align*}
We compute
\begin{align*}
	f_\pm^2g'_\pm & = \frac{1}{\sqrt{a^2\cos^2+b^2}}\Bigl[-(2a^2\cos^2 + b^2)a^2\cos^2 - 2a^2\cos^2(a^2\cos^2+b^2)\Bigr]\\
	&\quad \pm(2a^2\cos^2 + b^2)a\cos \pm 2a^3\cos^3.
\end{align*}
Now, in view of~\eqref{eq:periodic-nodoid}, the (orientation dependent) volume $\mathcal V$ of the rotational symmetric surface corresponding to one period can be computed as
\begin{align*}
	\mathcal V & = \pi\int_{-\frac{\pi}{2}}^\frac{\pi}{2}f_-^2g_-'\,\mathrm ds + \pi\int_{-\frac{\pi}{2}}^\frac{\pi}{2}f_+^2g_+'\,\mathrm ds = \pi\int_{-\pi}^\pi f_-^2g_-'\,\mathrm ds\\
	&=-4\pi a^2\int_{0}^\frac{\pi}{2} \frac{4a^2\cos^4(s) + 3b^2\cos^2(s)}{\sqrt{a^2\cos^2(s) + b^2}}\,\mathrm ds.
\end{align*}
This integral can be expressed in terms of complete elliptic integrals. For the purpose of this paper it is enough to notice that
\begin{equation}\label{eq:convergence_volume-nodoid}
	\mathcal V\xrightarrow{b\downarrow0}-\frac{4\pi}{3}(2a)^3.
\end{equation}
	
\section{Embedded Delaunay tori} \label{sec:d-tori}
\begin{figure}[h]
	\begin{minipage}[b]{0.45\linewidth}
		\centering
		\includegraphics[width=\textwidth]{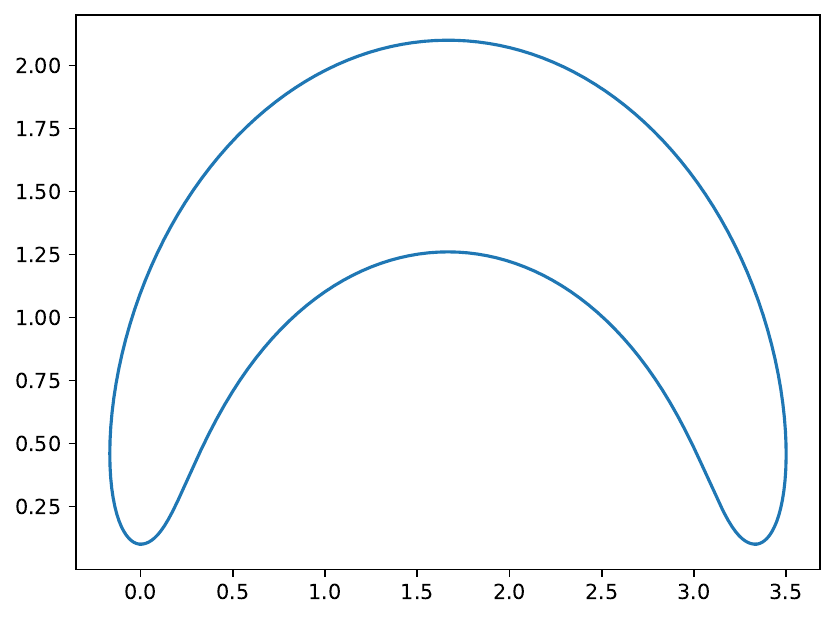}
		\caption{Profile curve Delaunay torus with $c=1.1$ and the bottom line being the axis of rotation.}
		\label{fig:d-torus}
	\end{minipage}
	\hspace{0.5cm}
	\begin{minipage}[b]{0.45\linewidth}
		\centering
		\includegraphics[width=\textwidth]{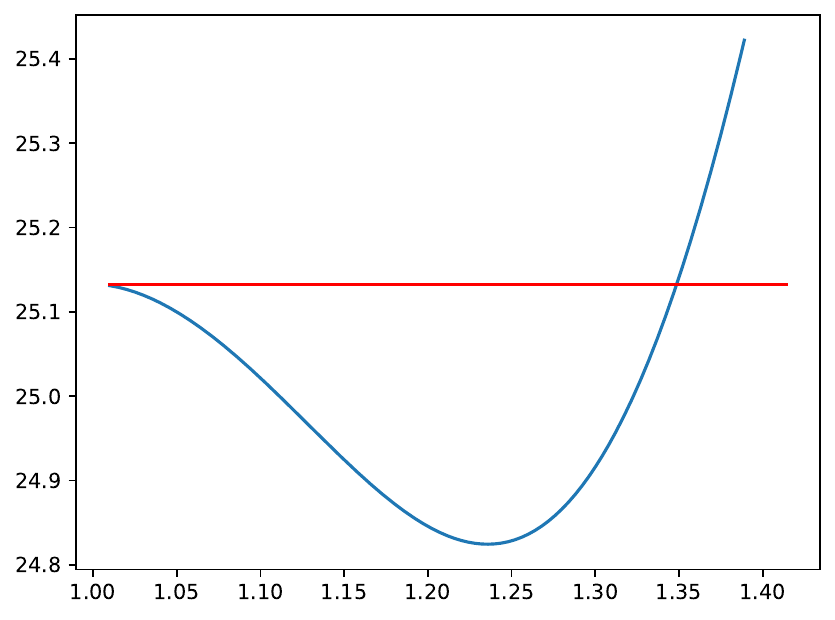}
		\caption{Energy curve for the family of Delaunay tori and $8\pi$ bound. \\}
		\label{fig:total_energy}
	\end{minipage}
\end{figure}
In this Section, we will construct the family of embedded Delaunay tori $\mathbb T_{\mathrm D,c}$ with $1<c<\gamma_0$ for some constant~$\gamma_0$ and we will prove Theorem~\ref{thm:8pi-bound}. Moreover, at the end of this section, we will give a proof of Corollaries~\ref{cor:existence_iso_constrained_tori} and~\ref{cor:conf-constrained}. 

Each Delaunay torus is a rotationally symmetric surface whose profile curve (see Figure~\ref{fig:d-torus}) consists of one period of an unduloid roulette (see Figure~\ref{fig:unduloid-2periods}) and one period of a nodoid roulette (see Figure~\ref{fig:nodoid-2periods}). The construction works as follows. Start with a one parameter family of patched nodoids (see Section~\ref{sec:nodoids}) running for one period and starting at the point where $f_+$ attains its minimum
according to \eqref{eq:nodoid-extrem_values}, where $a=1$ and $c>1$ is the free parameter, thus $b$ is given by $b = \sqrt{c^2 - 1}$ and the minimum according to \eqref{eq:nodoid-extrem_values} is $c-1$. Next, depending on the parameter $c$, find $a>y>0$ such that the unduloid corresponding to the ellipse with foci $(\pm y,0)$, width $2a$, and height $2b$ for $b = \sqrt{a^2 - y^2}$ (see Section~\ref{sec:unduloids}) running for one period and starting at 
the point where $f$ attains its minimum
$a - y$ according \eqref{eq:unduloid-extrem_values}, fits right into the given nodoid. That means the two end points where the patched nodoids reach their minimum need to match the two end points where the unduloid reaches its minimum. Notice that, in this way, the profile curve is $C^{1,1}$-regular. 
Indeed, choosing graphical representations around the glueing points (such representations exist since by \eqref{eq:unduloid-g_component}, $g'(0) = g'(2\pi)\neq 0$ and by \eqref{eq:nodoid-g_component}, $g_-'(0)\neq0$), it is clear that the glued curve is $C^1$-regular and, since both unduloids and nodoids induce constant mean curvature surfaces, their second derivative is naturally bounded near the gluing point leading to $C^{1,1}$-regularity, cf.\ \cite[Equation (3.3)]{MR1807955}.
The coordinates of the two patching points can be determined using Equations \eqref{eq:unduloid-extrem_values}, \eqref{eq:unduloid-extrinsic_length} for the unduloid and \eqref{eq:nodoid-extrem_values}, \eqref{eq:nodoid-extrinsic_length} for the nodoid. Thus, $a,y$ are given as the solution of the system of equations
\begin{empheq}[left=\empheqlbrace]{align} \label{eq:d-tori_balncing_cond1}
	4aE\Bigl(\frac{y}{a}\Bigr) & = 4c\left[E\Bigl(\frac{1}{c}\Bigr) - \Bigl(1 - 
	\frac{1}{c^2}\Bigr)K\Bigl(\frac{1}{c}\Bigr)\right]\\ \label{eq:d-tori_balncing_cond2}
	a - y & = c - 1.
\end{empheq}
Abbreviating 
\begin{equation*}
	 \varepsilon := c - 1, \qquad \Lambda := c\left[E\Bigl(\frac{1}{c}\Bigr) - \varepsilon \frac{c+1}{c^2}K\Bigl(\frac{1}{c}\Bigr)\right],
\end{equation*}
the system of equations reads as 
\begin{empheq}[left=\empheqlbrace]{align*} 
	(y + \varepsilon)E\Bigl(\frac{y}{y + \varepsilon}\Bigr) & = \Lambda\\
	 a & = y + \varepsilon.
\end{empheq}
Define the $C^1$-function
\begin{equation*}
	F:(1,\infty)\times(0,\infty)\to \mathbb R, \qquad F(c,y) = (y + \varepsilon(c))E\Bigl(\frac{y}{y + \varepsilon(c)}\Bigr) - \Lambda(c).
\end{equation*}
There holds 
\begin{equation*}
	\varepsilon \log\frac{4}{\sqrt{1 - 1/c^2}} = \varepsilon \log \frac{4c}{\sqrt{\varepsilon}\sqrt{c+1}} \leq \frac{4c}{\sqrt{c+1}}\sqrt{\varepsilon} \xlongrightarrow{c \downarrow 1} 0.	
\end{equation*}
Hence, by \eqref{eq:eInts-growth_of_K} and since $E(1)=1$, we have
\begin{equation} \label{eq:d-tori-limit_L}
	\lim_{c\to 1+}\Lambda(c) = 1.
\end{equation}
Moreover, 
\begin{align*}
\partial_c\Lambda & = [E(\textstyle\frac{1}{c}) - (1 - \frac{1}{c^2})K(\frac{1}{c})] + c[E(\frac{1}{c}) - K(\frac{1}{c})]c\frac{-1}{c^2}\\  
&\quad +  c\left[\textstyle\frac{-2}{c^3}K(\frac{1}{c}) - (1 - \frac{1}{c^2})[E(\frac{1}{c})/(1 - \frac{1}{c^2}) - K(\frac{1}{c})]c(\frac{-1}{c^2})\right] = E(\textstyle\frac{1}{c}) - K(\frac{1}{c})
\end{align*}
and, for $k = y/(y+\varepsilon)$ 
\begin{equation*}
	\partial_c(F+\Lambda) = E(k) + (y+\varepsilon)[E(k) - K(k)]\frac{y+\varepsilon}{y}\frac{-y}{(y+\varepsilon)^2} = K(k)
\end{equation*}
which implies
\begin{equation*}
	\partial_cF = K\Bigl(\frac{y}{y+\varepsilon}\Bigr) + K\Bigl(\frac{1}{c}\Bigr) - E\Bigl(\frac{1}{c}\Bigr).
\end{equation*}
Writing $k = y/(y+\varepsilon)$, there holds
\begin{equation*}
	\partial_yF=E(k) + (y + \varepsilon)[E(k) - K(k)]\frac{y+\varepsilon}{y}\frac{\varepsilon}{(y + \varepsilon)^2} = \Bigl(1 + \frac{\varepsilon}{y}\Bigr)E(k) -\frac{\varepsilon}{y}K(k).
\end{equation*}
Moreover, for fixed $y>0$,
\begin{equation}\label{eq:d-tori-limit_eps_times_K}
	\varepsilon \log \frac{4}{\sqrt{1 - y^2/(y + \varepsilon)^2}} = \varepsilon\log \frac{4(y + \varepsilon)}{\sqrt{2\varepsilon y + \varepsilon^2}} \leq \frac{4(y + \varepsilon)}{\sqrt{2y + \varepsilon}}\sqrt{\varepsilon}\xlongrightarrow{c \downarrow 1} 0
\end{equation} 
and thus, by~\eqref{eq:eInts-growth_of_K},
\begin{equation}\label{eq:part_der_F}
	\lim_{c \to 1+}\partial_yF(c,y)=1.
\end{equation}
In view of~\eqref{eq:d-tori-limit_L}, we can choose $1<\gamma_0<1+(2\pi)^{-1}$ such that $\Lambda(c)>1/2$ for all $1<c<\gamma_0$. For such $c$ it follows that $\varepsilon(c)\leq(2\pi)^{-1}$ and by~\eqref{eq:eInts-boundedness_E},
\begin{equation*}
	(y + \varepsilon(c))E\Bigl(\frac{y}{y+\varepsilon(c)}\Bigr)\leq\frac{1}{2},\qquad F(c,y)<0
\end{equation*} 
for all $y\leq(2\pi)^{-1}$. On the other hand, using ~\eqref{eq:eInts-boundedness_E}, we have $\Lambda(c)\leq \gamma_0\pi/2<3$ and thus $F(c,y)>0$ for all $y\geq3$. Therefore, by the intermediate value theorem, we see that for all $1<c<\gamma_0$ there exists 
\begin{equation}\label{eq:bound_on_y}
	\frac{1}{2\pi}\leq y(c)\leq 3
\end{equation}
such that $F(c,y(c)) = 0$.
Moreover, since $y(c)\geq(2\pi)^{-1}$, we can use~\eqref{eq:part_der_F} to see that for $c$ close to 1, $\partial_yF(c,y(c)) \neq 0$.
Therefore, we can apply the implicit function theorem to obtain $\gamma_0>1$ such that $y(c)$ is a smooth function in $c$ for $1<c<\gamma_0$ and
\begin{equation} \label{eq:d-tori-derivative_y}
	y' = \frac{E(\frac{1}{c}) - K(\frac{1}{c}) - K(\frac{y}{y + \varepsilon})}{(y + \varepsilon) E(\frac{y}{y + \varepsilon})-\varepsilon K(\frac{y}{y + \varepsilon})} \,y.
\end{equation}
Using \eqref{eq:unduloid-mean_curvature}, \eqref{eq:unduloid-area} for the unduloid and \eqref{eq:nodoid-mean_curvature}, \eqref{eq:nodoid-area} for the nodoid (with $c = y$,	$a = y + \varepsilon$), one obtains the Willmore energy of the Delaunay tori $\mathbb T_{\mathrm D,c}$:
\begin{equation} \label{eq:d-tori-energy_decomposition}
	\mathcal W(\mathbb T_{\mathrm D,c}) = \mathcal W_\mathrm{nod} + \mathcal W_\mathrm{und}
\end{equation}
where
\begin{equation*}
	\mathcal W_\mathrm{nod} = 2\pi(1+c)E\Bigl(\frac{2\sqrt{c}}{1+c}\Bigr), \qquad \mathcal W_\mathrm{und} = 2\pi\Bigl(1 + \frac{y}{y + \varepsilon}\Bigr)E\Bigl(\frac{2\sqrt{y(y + \varepsilon)}}{2y + \varepsilon}\Bigr).
\end{equation*}
Using the uniform bound on $y$ in~\eqref{eq:bound_on_y}, we infer 
by~\eqref{eq:d-tori-limit_eps_times_K}, \eqref{eq:eInts-growth_of_K}, and $E(1)=1$ that
\begin{equation*}
	\lim_{c\to1+}\frac{y}{y + \varepsilon} = 1,\qquad \lim_{c\to1+}y = 1, \qquad \lim_{c \to 1+}\varepsilon K\Bigl(\frac{y}{y+\varepsilon}\Bigr)=0, \qquad \lim_{c\to1+} \mathcal W(\mathbb T_{\mathrm D,c}) = 8\pi.
\end{equation*}

Next, we show that $\partial_c \mathcal W(\mathbb T_{\mathrm D,c}) < 0$ for $c$ close to 1 which then implies $\mathcal W(\mathbb T_{\mathrm D,c}) < 8\pi$ for $c$ close to 1. First, we compute $\partial_c \mathcal W_\mathrm{nod}$. For this purpose, let $k = 2\sqrt{c}/(1 + c)$. Then,
\begin{equation*}
	\partial_ck= \frac{1}{(1+c)\sqrt{c}} - \frac{2\sqrt{c}}{(1+c)^2} = \frac{1 - c}{(1+c)^2\sqrt{c}}
\end{equation*}
and 
\begin{align*}
	\partial_c \mathcal W_\mathrm{nod} & = 2\pi E(k) + 2\pi(1+c)[E(k) - K(k)]\frac{1+c}{2\sqrt{c}}\frac{1 - c}{(1+c)^2\sqrt{c}}\\
	& = \pi\Bigl(\Bigl(1+\frac{1}{c}\Bigr)E(k) + \Bigl(1 - \frac{1}{c}\Bigr)K(k)\Bigr).
\end{align*}
By the Gau{\ss} transformation~\eqref{eq:eInts-Gauss_trans} there hold
\begin{equation*}
		k' = \sqrt{1 - \frac{4c}{(1+c)^2}} = \frac{c - 1}{c+1}, \qquad k_1 = \frac{1 - k'}{1 + k'} = \frac{1}{c},
\end{equation*}
and
\begin{align*}
	\textstyle\frac{\partial_c \mathcal W_\mathrm{nod}}{\pi} = (1 + \frac{1}{c})\left[(1 + \frac{c - 1}{c+1})E(\frac{1}{c}) - \frac{c-1}{c+1}(1+ \frac{1}{c})K(\frac{1}{c})\right] + (1-\frac{1}{c})(1 + \frac{1}{c})K(\frac{1}{c}) = 2E(\frac{1}{c}).	
\end{align*}
Hence,
\begin{equation} \label{eq:d-tori-der_energy_nod}
	\partial_c\mathcal W_\mathrm{nod} = 2\pi E\Bigl(\frac{1}{c}\Bigr).
\end{equation}

In order to compute $\partial_y \mathcal W_\mathrm{und}$, let $k = 2\sqrt{y(y + \varepsilon)}/(2y + \varepsilon)$. Then, there hold
\begin{equation*}
	\textstyle\partial_y k = \frac{2y + \varepsilon}{(2y + \varepsilon)\sqrt{y(y + \varepsilon)}} - \frac{4\sqrt{y(y + \varepsilon)}}{(2y + \varepsilon)^2} = \frac{\varepsilon^2}{(2y + \varepsilon)^2\sqrt{y(y + \varepsilon)}}, \qquad \partial_y \frac{y}{y + \varepsilon} = \frac{\varepsilon}{(y + \varepsilon)^2}
\end{equation*}
and
\begin{align*}
	\partial_y \mathcal W_\mathrm{und} & = 2\pi \textstyle\frac{\varepsilon}{(y + \varepsilon)^2}E(k) + 2\pi \frac{2y+\varepsilon}{y + \varepsilon}[E(k) - K(k)] \frac{2y + \varepsilon}{2\sqrt{y(y + \varepsilon)}}\frac{\varepsilon^2}{\sqrt{y(y + \varepsilon)}(2y + \varepsilon)^2}\\
	&=\textstyle \frac{\pi\varepsilon}{(y + \varepsilon)^2}\bigl(2E(k) + \frac{\varepsilon}{y}[E(k) - K(k)]\bigr) = \frac{\pi\varepsilon}{y(y + \varepsilon)^2}\bigl((2y + \varepsilon)E(k) - \varepsilon K(k)\bigr).
\end{align*}
By the Gau{\ss} transformation~\eqref{eq:eInts-Gauss_trans} there hold 
\begin{equation*}
	k' = \sqrt{1 - \frac{4y(y + \varepsilon)}{(2y + \varepsilon)^2}} = \frac{\varepsilon}{2y + \varepsilon}, \qquad k_1 = \frac{1 - k'}{1 + k'} = \frac{y}{y + \varepsilon}
\end{equation*}
and 
\begin{equation} \label{eq:d-tori-Gauss_trans_energy_und}
	\begin{split}
		&(2y + \varepsilon) E(k) - \varepsilon K(k) = (2y + \varepsilon)[\textstyle (1 + \frac{\varepsilon}{2y + \varepsilon})E(\frac{y}{y + \varepsilon}) - \frac{\varepsilon}{2y + \varepsilon}(1 + \frac{y}{y + \varepsilon})K(\frac{y}{y + \varepsilon})] \\
		&\qquad \textstyle - \varepsilon(1 + \frac{y}{y + \varepsilon})K(\frac{y}{y + \varepsilon}) = \textstyle 2\bigl((y + \varepsilon)E(\frac{y}{y + \varepsilon}) - \varepsilon(1 + \frac{y}{y + \varepsilon})K(\frac{y}{y + \varepsilon})\bigr).
	\end{split}
\end{equation}
Therefore,
\begin{equation} \label{eq:d-tori-der_energy_und_y}
	\partial_y \mathcal W_\mathrm{und} = \frac{2\pi\varepsilon}{y(y + \varepsilon)^2}\left[(y + \varepsilon)E\Bigl(\frac{y}{y + \varepsilon}\Bigr) - \varepsilon\Bigl(1 + \frac{y}{y + \varepsilon}\Bigr)K\Bigl(\frac{y}{y + \varepsilon}\Bigr)\right].
\end{equation}
Abbreviate $z = y'/y$, and $k = y/(y+\varepsilon)$. Then,  \eqref{eq:d-tori-der_energy_und_y} and \eqref{eq:d-tori-derivative_y} imply  
\begin{equation} \label{eq:d-tori-der_energy_und_y_times_y'}
	\begin{split} 
		\partial_y \mathcal{W}_\mathrm{und}\cdot y' &= \textstyle \frac{2\pi\varepsilon}{y(y + \varepsilon)^2}\left[(y + \varepsilon)E(k) - \varepsilon(1 + \frac{y}{y + \varepsilon})K(k)\right]\cdot \frac{E(\frac{1}{c}) - K(\frac{1}{c}) - K(k)}{(y + \varepsilon) E(k)-\varepsilon K(k)} \,y\\
		&= \textstyle \frac{2\pi\varepsilon}{(y + \varepsilon)^2}\left[E(\frac{1}{c}) - K(\frac{1}{c}) - K(\frac{y}{y + \varepsilon})\right] - \frac{2\pi\varepsilon}{(y + \varepsilon)^2}\frac{\varepsilon y}{y+\varepsilon} z K(\frac{y}{y + \varepsilon}).
	\end{split}
\end{equation}

Finally we compute $\partial_\varepsilon \mathcal W_\mathrm{und}$. For this purpose let $k = 2\sqrt{y(y + \varepsilon)}/(2y + \varepsilon)$. Then, there hold
\begin{equation*}
	\partial_\varepsilon k = \frac{y}{(2y + \varepsilon)\sqrt{y(y + \varepsilon)}} - \frac{2\sqrt{y(y + \varepsilon)}}{(2y + \varepsilon)^2} = \frac{-y\varepsilon}{(2y + \varepsilon)^2\sqrt{y(y + \varepsilon)}}
\end{equation*}
and 
\begin{align*}
	\partial_\varepsilon \mathcal W_\mathrm{und} & = \textstyle2\pi \frac{-y}{(y + \varepsilon)^2} E(k) + 2\pi\frac{2y + \varepsilon}{y + \varepsilon}[E(k) -K(k)]\frac{2y + \varepsilon}{2\sqrt{y(y + \varepsilon)}}\frac{-y\varepsilon}{(2y + \varepsilon)^2\sqrt{y(y + \varepsilon)}} \\
	& = \textstyle-\frac{\pi}{(y + \varepsilon)^2}\bigl((2y + \varepsilon)E(\frac{y}{y + \varepsilon}) - \varepsilon K(\frac{y}{y + \varepsilon})\bigr).
\end{align*}
Thus, by \eqref{eq:d-tori-Gauss_trans_energy_und},
\begin{equation*}
	\partial_\varepsilon \mathcal W_\mathrm{und} = - \frac{2\pi}{(y + \varepsilon)^2}\left[(y + \varepsilon)E\Bigl(\frac{y}{y + \varepsilon}\Bigr) - \varepsilon\Bigl(1 + \frac{y}{y + \varepsilon}\Bigr)K\Bigl(\frac{y}{y + \varepsilon}\Bigr)\right].
\end{equation*}
Recall that, by the choice of $y$, there holds 
\begin{equation*}
	(y + \varepsilon)E\Bigl(\frac{y}{y + \varepsilon}\Bigr) = c\left[E\Bigl(\frac{1}{c}\Bigr) - \Bigl(1 - \frac{1}{c^2}\Bigr)K\Bigl(\frac{1}{c}\Bigr)\right].
\end{equation*}
Therefore, 
\begin{equation} \label{eq:d-tori-der_energy_und_eps}
	\partial_\varepsilon \mathcal W_\mathrm{und} = - \frac{2\pi}{(y + \varepsilon)^2}\left[cE\Bigl(\frac{1}{c}\Bigr) -\varepsilon\frac{c+1}{c}K\Bigl(\frac{1}{c}\Bigr)- \varepsilon\Bigl(1 + \frac{y}{y + \varepsilon}\Bigr)K\Bigl(\frac{y}{y + \varepsilon}\Bigr)\right].
\end{equation}

Putting \eqref{eq:d-tori-der_energy_nod}, \eqref{eq:d-tori-der_energy_und_y_times_y'}, and \eqref{eq:d-tori-der_energy_und_eps} into \eqref{eq:d-tori-energy_decomposition}, it follows
\begin{align*}
	& \partial_c \mathcal W(\mathbb T_{\mathrm D,c}) = \partial_c \mathcal W_\mathrm{nod} + \partial_\varepsilon \mathcal W_\mathrm{und} + \partial_y \mathcal W_\mathrm{und} \cdot y'\\
	&\quad = 2\pi \textstyle E(\frac{1}{c}) - \frac{2\pi}{(y + \varepsilon)^2}\left[cE\Bigl(\frac{1}{c}\Bigr) -\varepsilon\frac{c+1}{c}K\Bigl(\frac{1}{c}\Bigr)- \varepsilon\Bigl(1 + \frac{y}{y + \varepsilon}\Bigr)K\Bigl(\frac{y}{y + \varepsilon}\Bigr)\right] \\
	&\qquad + \textstyle \frac{2\pi\varepsilon}{(y + \varepsilon)^2}\left[E(\frac{1}{c}) - K(\frac{1}{c}) - K(\frac{y}{y + \varepsilon})\right] - \frac{2\pi\varepsilon^2y}{(y + \varepsilon)^3} z K(\frac{y}{y + \varepsilon}) \\
	& \quad = 2\pi \textstyle E(\frac{1}{c})\left[1 - \frac{c}{(y + \varepsilon)^2} + \frac{\varepsilon}{(y + \varepsilon)^2}\right] + \frac{2\pi\varepsilon}{(y + \varepsilon)^2}K(\frac{1}{c})\left[\frac{c+1}{c} - 1\right] + \frac{2\pi\varepsilon y}{(y + \varepsilon)^3}K(\frac{y}{y + \varepsilon})[1-\varepsilon z] \\
	& \quad = \textstyle 2\pi (1 - \frac{1}{(y + \varepsilon)^2}) E(\frac{1}{c}) + \frac{2\pi}{c(y + \varepsilon)^2}\varepsilon K(\frac{1}{c}) + \frac{2\pi y}{(y + \varepsilon)^3}(1 - \varepsilon z)\varepsilon K(\frac{y}{y + \varepsilon}).
\end{align*}
Therefore,
\begin{equation} \label{eq:d-tori-der_total_energy}
	\frac{\partial_c \mathcal W(\mathbb T_{\mathrm D,c})}{2\pi} = a_1 \Bigl(1 - \frac{1}{(y + \varepsilon)^2}\Bigr) + a_2 \varepsilon K\Bigl(\frac{1}{c}\Bigr) + a_3 \varepsilon K\Bigl(\frac{y}{y + \varepsilon}\Bigr)
\end{equation}
with $a_1,a_2,a_3$ depending on $c$. Recalling that
\begin{equation}\label{eq:limits}
	\textstyle \varepsilon K(\frac{1}{c}) \xlongrightarrow{c\downarrow 1} 0,\quad \varepsilon K(\frac{y}{y + \varepsilon})\xlongrightarrow{c\downarrow 1} 0, \quad \varepsilon y'\xlongrightarrow{c\downarrow 1}0, \quad y'\xlongrightarrow{c\downarrow 1}-\infty, \quad y\xlongrightarrow{c\downarrow 1}1
\end{equation}
we in particular see $a_1,a_2,a_3 \to 1$ as $c\to1$ and hence $\partial_c \mathcal W(\mathbb T_{\mathrm D,c})\to 0$ as $c\to1$. We claim that
\begin{equation}\label{eq:d-tori-l'hospital_limit}
	\lim_{c \to 1+} \frac{\varepsilon K(\frac{1}{c}) + \varepsilon K(\frac{y}{y + \varepsilon})}{1 - \frac{1}{(y + \varepsilon)^2}} = -\frac{1}{2}.
\end{equation}
In view of \eqref{eq:limits}, it follows
\begin{gather*}
	\textstyle\frac{\varepsilon}{1 + y'}\partial_c K(\frac{1}{c}) = \frac{\varepsilon}{1 + y'}\left[\frac{E(\frac{1}{c})}{1 - \frac{1}{c^2}} - K(\frac{1}{c})\right]c\frac{-1}{c^2} \xlongrightarrow{c\downarrow 1}0,\\
	\textstyle\partial_c \frac{y}{y+\varepsilon} = y'\partial_y \frac{y}{y+\varepsilon} + \partial_\varepsilon \frac{y}{y+\varepsilon} = \frac{\varepsilon y' - y}{(y + \varepsilon)^2}\xlongrightarrow{c\downarrow 1} - 1,\\
	\textstyle\frac{\varepsilon}{1 + y'}\partial_c K(\frac{y}{y+\varepsilon}) = \frac{\varepsilon}{1 + y'}\left[\frac{E(\frac{y}{y + \varepsilon})}{1 - \frac{y^2}{(y + \varepsilon)^2}} - K(\frac{y}{y + \varepsilon})\right]\frac{y + \varepsilon}{y}\partial_c \frac{y}{y+\varepsilon}\xlongrightarrow{c\downarrow 1}0.
\end{gather*}
Thus, by L'H\^ospital's rule and~\eqref{eq:d-tori-derivative_y},
\begin{equation*}
	 \lim_{c \to 1+} \frac{\varepsilon K(\frac{1}{c}) + \varepsilon K(\frac{y}{y + \varepsilon})}{1 - \frac{1}{(y + \varepsilon)^2}} = \lim_{c \to 1+} \frac{K(\frac{1}{c}) +  K(\frac{y}{y + \varepsilon})}{\frac{2}{(y + \varepsilon)^3}[1 + y']} = -\frac{1}{2}
\end{equation*}
which proves \eqref{eq:d-tori-l'hospital_limit}. By \eqref{eq:d-tori-l'hospital_limit} and \eqref{eq:d-tori-der_total_energy}, one infers
\begin{equation*}
	\partial_c \mathcal W(\mathbb T_{\mathrm{D},c}) < 0 \qquad \text{for $c$ close to 1}.
\end{equation*}
Therefore, for some $\gamma_0 > 1$, there holds
\begin{equation}\label{eq:d-tori-8pi_bound}
	\mathcal W(\mathbb T_{\mathrm D,c}) < 8\pi \qquad \text{whenever } 1< c<\gamma_0.
\end{equation}

By the Li--Yau inequality~\cite{MR674407} we see that the Delaunay tori $\mathbb T_{\mathrm{D},c}$ are embedded, the profile curve is a simple closed curve, the genus is $1$, and the algebraic volume~\eqref{eq:algebraic_volume} is the enclosed volume (up to orientation). Combining \eqref{eq:unduloid-area}, \eqref{eq:volume_convergence-unduloid}, \eqref{eq:nodoid-area}, and \eqref{eq:convergence_volume-nodoid}, we infer
\begin{equation*}
	\lim_{c \to 1+} \iso(\mathbb T_{\mathrm{D},c}) = \infty.
\end{equation*}
Together with \eqref{eq:d-tori-8pi_bound}, this proves Theorem~\ref{thm:8pi-bound}. It remains to mention that, 
since the profile curve is $C^{1,1}$-regular as discussed at the beginning of this section, the resulting surfaces of revolution $\mathbb T_{\mathrm{D},c}$ are also $C^{1,1}$-regular. \qed\medskip 

\textbf{Proof of Corollary~\ref{cor:existence_iso_constrained_tori}.} 
Let $\mathbb T^2$ be an abstract 2-dimensional torus. Denote with $\mathcal E_{\mathbb S^2}$ and $\mathcal E_{\mathbb T^2}$
the spaces of Lipschitz immersions of $\mathbb S^2$ and $\mathbb T^2$ into~$\mathbb R^3$ as defined in Section~2.2 of~\cite{MR3176354}.
By \textsc{Schygulla} \cite{Schygulla} and \cite[Theorem~1.1]{MR3176354} the following holds true. For each $\sigma>\sqrt[3]{36\pi}$, there exists a smoothly embedded spherical surface $\mathbb S_{\mathrm S,\sigma}$ with
\begin{equation*}
	\mathcal W(\mathbb S_{\mathrm S,\sigma}) = \beta_0(\sigma) := \inf\{\mathcal W(\vec \Phi) \mid \vec{\Phi} \in \mathcal E_{\mathbb S^2},\,\iso(\vec\Phi) = \sigma\}.
\end{equation*}
Moreover, by Theorem 1.6 in \textsc{Keller--Mondino--Rivi\`ere}~\cite{MR3176354}, for each $\sigma_0$ in the set
\begin{equation*}
	I_1:=\left\{\sigma \in \mathbb R \,\middle|\, \inf_{\substack{\vec \Phi \in \mathcal E_{\mathbb T^2} \\ \iso(\vec \Phi) = \sigma}} \mathcal W(\vec \Phi) < \min\{8\pi, 2\pi^2 + \beta_0(\sigma) - 4\pi\}\right\} \subset (\sqrt[3]{36\pi},\infty)
\end{equation*}
there exists a smoothly embedded torus $\Sigma_0$ in $\mathbb R^3$ with
\begin{equation*}
	\mathcal W(\Sigma_0) = \beta_1(\sigma_0):= \inf_{\substack{\vec \Phi \in \mathcal E_{\mathbb T^2} \\ \iso(\vec \Phi) = \sigma_0}} \mathcal W(\vec \Phi).
\end{equation*}
From Corollary 1.5 in \textsc{Mondino--Scharrer}~\cite{mondino2020strict} it follows
\begin{equation*}
	I_1 = \{\sigma \in \mathbb R \mid \beta_1(\sigma) < 8\pi\}.
\end{equation*}
Recall that the function $\beta_1(\cdot)$ is non-decreasing on the set $I_1$. This can be shown using M\"obius transformations, see for instance \cite[Proposition~1]{scharrer:PhD} and  \cite[Theorem~3.1]{yu2020uniqueness}. Therefore, the set $I_1$ is in fact an interval. Moreover, each $C^{1,1}$-embedding of $\mathbb T^2$ into~$\mathbb R^3$ is in particular a Lipschitz immersion, i.e.\ a member of the space $\mathcal E_{\mathbb T^2}$. Thus, by Theorem~\ref{thm:8pi-bound},
\begin{equation*}
	I_1 = (\sqrt[3]{36\pi},\infty)
\end{equation*}
which proves Corollary~\ref{cor:existence_iso_constrained_tori}.\qed\medskip 

\textbf{Proof of Corollary~\ref{cor:conf-constrained}.} 
Denote with $\alpha_c:\mathbb S^2 \to \mathbb R^2$ the $C^{1,1}$-regular profile curve of the Delaunay tori $\mathbb T_{\mathrm D,c}$ corresponding to $1<c<\gamma_0$. In view of~\eqref{eq:unduloid-extrem_values} and~\eqref{eq:nodoid-extrem_values}, we have indeed that $\alpha_c$ takes values in the half space $\mathbb H^2:=\{(x,y)\in\mathbb R^2\mid x>0\}$. Let~$\alpha_c^\varepsilon$ be the mollification of $\alpha_c$ for $\varepsilon>0$ and let $T_{\alpha_c^\varepsilon}$ be the corresponding rotationally symmetric torus. Then, 
\begin{equation*}
	T_{\alpha_c^\varepsilon} \to \mathbb T_{\mathrm D,c}\qquad \text{as $\varepsilon \downarrow0$ in $W^{2,p}$ for all $p>1$}.
\end{equation*}
In particular, for $\varepsilon >0$ small, we have that $\alpha_c^\varepsilon\in C^{\infty}(\mathbb S^1,\mathbb H^2)$ and, by Theorem~\ref{thm:8pi-bound}, $\mathcal W(T_{\alpha_c^\varepsilon})<8\pi$. Following \cite[Equation~(2.6)]{dall2020willmore}, we define 
\begin{equation*}
	\mathcal L_{\mathbb H^2}(\gamma) := \int_0^1\frac{|\gamma'|}{\gamma^1}\mathrm dx\qquad \text{for all $\gamma=(\gamma^1,\gamma^2)\in C^1(\mathbb S^1,\mathbb H^2)$.} 
\end{equation*}
Using \eqref{eq:unduloid-f_component} and \eqref{eq:speed-unduloid}, we see
\begin{equation}\label{eq:conf_class-unduloid}
\int_0^{2\pi}\frac{|(f,g)'|}{f}\,\mathrm dt = a\int_{0}^{2\pi}\frac{\mathrm dt}{\sqrt{a^2-c^2\cos^2(t)}}\xrightarrow{b\downarrow0}\int_0^{2\pi}\frac{\mathrm dt}{|\sin(t)|}=\infty.
\end{equation}
Similarly, using \eqref{eq:nodoid-f_component}, \eqref{eq:nodoid-speed}, and $\int\frac{1}{\sinh(t)} = \log |\tanh(t/2)| + C$, 
\begin{equation}\label{eq:conf_class-nodoid}
	\int_{-\infty}^\infty\frac{|(f_\pm,f_\pm)'|}{f_\pm}\,\mathrm dt = a\int_{-\infty}^\infty\frac{\mathrm dt}{\sqrt{c^2\cosh^2(t) - a^2}} \xrightarrow{b\downarrow0}\int_{-\infty}^\infty\frac{\mathrm dt}{|\sinh(t)|}=\infty.
\end{equation}
Hence, combining \eqref{eq:conf_class-unduloid} and \eqref{eq:conf_class-nodoid}, it follows that
\begin{equation*}
	\mathcal L_{\mathbb H^2}(\alpha_c) \to \infty \qquad \text{as $c\downarrow 1$.}
\end{equation*}
Moreover, for all $1<c<\gamma_0$, there holds
\begin{equation*}
	\mathcal L_{\mathbb H^2}(\alpha_c^\varepsilon) \to \mathcal L_{\mathbb H^2}(\alpha_c) \qquad \text{as $\varepsilon\downarrow 0$.}
\end{equation*}
Finally, since by \cite[Proposition 4.2.]{dall2020willmore}, the conformal class $\omega(T_\gamma)$ of any rotationally symmetric torus  $T_\gamma$ corresponding to a profile curve $\gamma\in C^\infty(\mathbb S^1,\mathbb H^2)$ is given by
\begin{equation*}
	\omega(T_\gamma) = \begin{cases}
		i\frac{1}{2\pi}\mathcal L_{\mathbb H^2}(\gamma) & \mathcal L_{\mathbb H^2}(\gamma)\geq 2\pi\\
		i\frac{2\pi}{\mathcal L_{\mathbb H^2}(\gamma)} & \mathcal L_{\mathbb H^2}(\gamma)<2\pi
	\end{cases},
\end{equation*}
the conclusion follows. \qed

\section{Delaunay spheres of high isoperimetric ratio} \label{sec:d-spheres}
\begin{figure}[h]
	\begin{minipage}{0.45\linewidth}
		\centering
		\includegraphics[width=\textwidth]{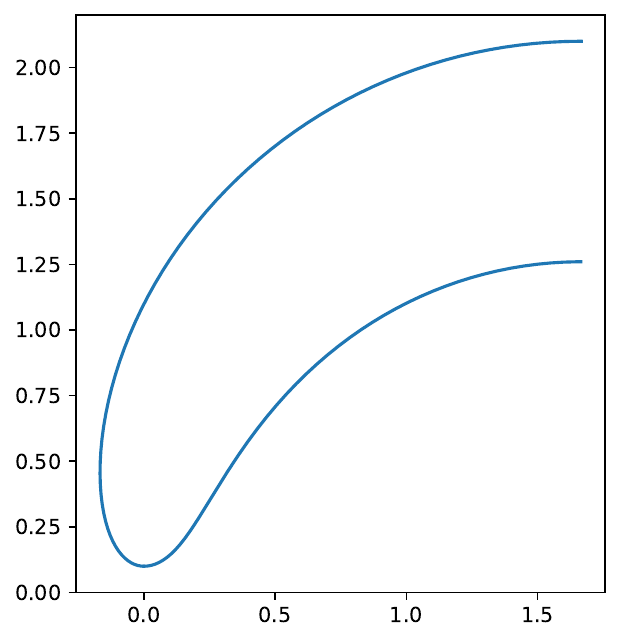}
		\caption{Profile curve of half a Delaunay torus with $c=1.1$}
		\label{fig:half_d-torus}
	\end{minipage}
	\hspace{0.5cm}
	\begin{minipage}{0.45\linewidth}
		\centering
		\includegraphics[width=\textwidth]{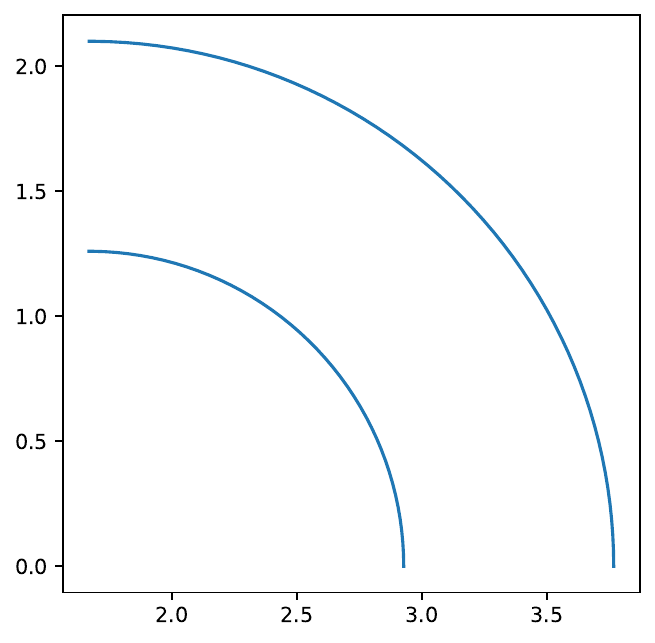}
		\caption{Concentric quarter circles fitting into half a Delaunay torus.}
		\label{fig:quarter}
	\end{minipage}
\end{figure}
\begin{figure}[h] 
	\centering
	\includegraphics[width=0.45\textwidth]{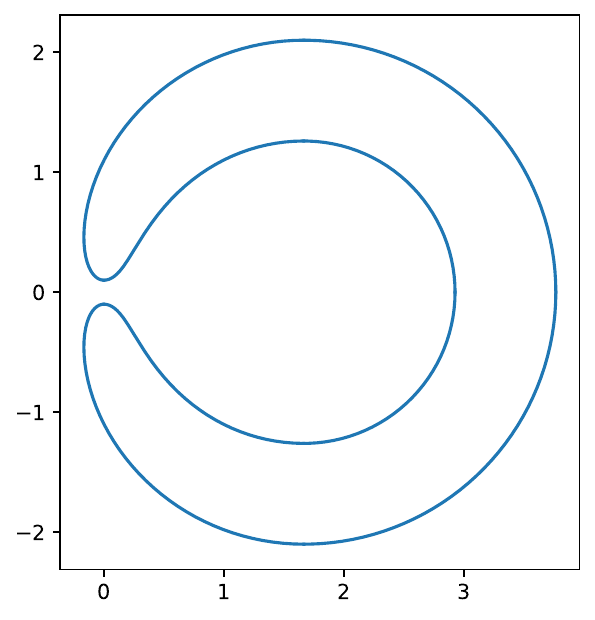}
	\caption{Profile of a Delaunay sphere with $c=1.1$}
	\label{fig:d-sphere}
\end{figure}

The Delaunay tori $\mathbb T_{\mathrm D,c}$ corresponding to $1<c<\gamma_0$ can be used to construct spheres with analogous properties. The first part of the construction works just like the construction of the Delaunay tori only that now, both the nodoid and the unduloid run only for half a period instead of one full period. To be more precise, both the nodoid and the unduloid now only run from 
the point where $f_+$, $f$ attain
their minimum according to \eqref{eq:unduloid-extrem_values}, \eqref{eq:nodoid-extrem_values} up to the point where $f_-$, $f$
reach their maximum (according to \eqref{eq:unduloid-extrem_values}, \eqref{eq:nodoid-extrem_values}) but not until $f_+$, $f$ reach their minimum again. This results in half a Delaunay torus, see Figure~\ref{fig:half_d-torus}. Notice that unduloids and nodoids are symmetric around their maxima ($t=\pi$ in \eqref{eq:unduloid-f_component}, \eqref{eq:unduloid-g_component} for unduloids; $t=0$ in \eqref{eq:nodoid-f_component}, \eqref{eq:nodoid-g_component} for nodoids). Thus, the Willmore energy of this particular half of a Delaunay torus is indeed half the Willmore energy of a whole Delaunay torus. Let $c,y,a$ be the balancing parameters according to \eqref{eq:d-tori_balncing_cond1} and \eqref{eq:d-tori_balncing_cond2}. Then the maxima of the nodoid and the unduloid are given by $c+1$ and $a+y$, respectively. Next, take two concentric circular sectors with radii $c+1$ and $a+y$ both of which being one quarter of a full circle, see Figure~\ref{fig:quarter}. Choose the centre of the two circular sectors at $L/2$ on the axis of rotation of the half Delaunay torus, where $L = 4aE(y/a)$ (see \eqref{eq:unduloid-extrinsic_length}). Then, the two circular sectors fit right into the half Delaunay torus, resulting in a $C^{1,1}$-curve. Since the two circular sectors meet the axis of rotation perpendicular, the resulting surface of revolution is $C^{1,1}$-regular too. It is of sphere type. The full profile can be seen in Figure~\ref{fig:d-sphere}. The resulting family of surfaces is called \emph{Delaunay spheres}. Since $a+y, c+1 \xlongrightarrow{c \downarrow 1} 2$, the Delaunay spheres converge as varifolds to a sphere of multiplicity 2 as $c \to 1$. Their Willmore energy is given by 
\begin{equation*}
\mathcal W(\mathbb S_{\mathrm D,c}) = 4\pi + \frac{\mathcal W(\mathbb T_{\mathrm D,c})}{2}.
\end{equation*}  
Similarly as in the proof of Theorem~\ref{thm:8pi-bound}, we can combine \eqref{eq:unduloid-area}, \eqref{eq:volume_convergence-unduloid}, \eqref{eq:nodoid-area}, and \eqref{eq:convergence_volume-nodoid} to see $\iso(\mathbb S_{\mathrm D,c}) \to \infty$ as $c\to1$ which concludes the proof of Theorem~\ref{thm:d-spheres}. \qed

\section{Helfrich tori with small spontaneous curvature} \label{sec:Helfrich_tori}
Suppose $A_0,V_0>0$ satisfy the isoperimetric inequality: $A_0^3 > 36\pi V_0^2$ and let 
\begin{equation*}
	\mathcal S_1(A_0,V_0) = \{f\in\mathcal S_1 \mid \area(f)=A_0,\,\vol(f)=V_0\}.
\end{equation*}
Notice that if $c_0 = 0$, then the Helfrich functional reduces to the Willmore functional: $\mathcal H_{0} = \mathcal W$. Moreover, one can show (see \cite[Equation~4.26]{MR4076069})
\begin{equation} \label{eq:Helfrich-continuity}
 	\left|\inf_{f\in\mathcal S_1(A_0,V_0)}\sqrt{\mathcal W(f)} - \inf_{f\in\mathcal S_1(A_0,V_0)}\sqrt{\mathcal H_{c_0}(f)} \right|\leq |c_0|\sqrt{A_0}.
\end{equation}
In particular, the minimal Helfrich energy is continuous with respect to $c_0$ at $c_0 = 0$. For the case $c_0 = 0$, existence of smoothly embedded minimisers with given fixed area and volume corresponds to Corollary 1.3 in~\cite{MR3176354}. Corollary~\ref{cor:Helfrich-tori} states that minimisers remain embedded for $c_0$ close to zero. Moreover, (by the choice of $\varepsilon(A_0,V_0)$ in \eqref{eq:Helfrich-epsilon}) a minimising sequence for $c_0=0$ has the same uniform bounds on the Willmore energy as a minimising sequence for $c_0$ close to zero. Indeed, we will see that the compactness proof in~\cite{MR3176354} still works for $c_0$ close to zero. However, for general $c_0$, minimisers are no longer embedded, see~\cite{MR4076069}. The following proof is a combination of four independent results: the two strict inequalities Theorem~\ref{thm:8pi-bound} and \cite[Corollary~1.5]{mondino2020strict} are needed to deduce that $\varepsilon(A_0,V_0)$ (as defined in \eqref{eq:Helfrich-epsilon}) is strictly positive; then, one can apply the compactness proof of~\cite{MR3176354}; finally, one can conclude the regularity from~\cite{MR4076069} (after \textsc{Rivi\`ere}~\cite{MR2430975}). 
 
In order to prove Corollary~\ref{cor:Helfrich-tori}, let 
$\vec \Phi_k$ be a minimising sequence of 
\begin{equation*}
	\inf\{\mathcal H_{c_0}(\vec \Phi) \mid \vec \Phi \in \mathcal E_{\mathbb T^2},\,\area(\vec \Phi) = A_0,\,\vol(\vec \Phi) = V_0\}.
\end{equation*}
Recall the definition of $\varepsilon(A_0,V_0)$ in Equation \eqref{eq:Helfrich-epsilon}:
\begin{equation*}
	\varepsilon(A_0,V_0):= \frac{\sqrt{\min\{8\pi, 2\pi^2 + \beta_0(A_0/V_0^{2/3}) - 4\pi\}} - \sqrt{\beta_1(A_0/V_0^{2/3})}}{2\sqrt{A_0}}.
\end{equation*}  
By Theorem~\ref{thm:8pi-bound} and \cite[Corollary~1.5]{mondino2020strict}, there holds $\varepsilon(A_0,V_0) > 0$. Using the continuity property \eqref{eq:Helfrich-continuity} one can show for $|c_0|< \varepsilon(A_0,V_0)$ that (see the proof of Lemma 4.4 in~\cite{MR4076069})
\begin{equation*}
	\limsup_{k\to\infty} \sqrt{\mathcal W(\vec \Phi_k)} < \sqrt{\beta_1(A_0/V_0^{2/3})} + 2\varepsilon(A_0,V_0)\sqrt{A_0}.
\end{equation*}
Thus, by the definition of $\varepsilon(A_0,V_0)$ it follows 
\begin{equation} \label{eq:Helfrich-8pi_bound}
	\limsup_{k\to\infty}\mathcal W(\vec \Phi_k) < 8\pi 
\end{equation}
and 
\begin{equation} \label{eq:Helfrich-strict_inequality}
	\limsup_{k\to\infty}\mathcal W(\vec \Phi_k) < 2\pi^2 + \beta_0(A_0/V_0^{2/3}) -4\pi. 
\end{equation}
In Section 4.3 of~\cite{MR3176354} it is shown that due to the strict inequality in \eqref{eq:Helfrich-strict_inequality}, the conformal factors of $\vec \Phi_k$ are bounded away from finitely many concentration points $a_1,\ldots,a_N$ in $\mathbb T^2$. Hence, by the uniform energy bound in \eqref{eq:Helfrich-8pi_bound}, there exists $\vec\Phi_\infty \in\mathcal E_{\mathbb T^2}$ such that (after passing to a subsequence and after re-parametrising) for all $\delta>0$
\begin{equation}\label{eq:Helfrich-weak_convergence}
	\vec \Phi_k \to \vec \Phi_\infty \qquad \text{as $k\to\infty$ weakly in $W^{2,2}(\mathbb T^2\setminus \textstyle\bigcup_{i=1}^NB_\delta(a_i),\mathbb R^3)$.}
\end{equation}
Moreover, in Section 4.2 of~\cite{MR3176354} it is shown that due to \eqref{eq:Helfrich-8pi_bound}, there hold that
\begin{equation*}
	\lim_{k\to\infty}\vol(\vec \Phi_k) = \vol(\vec \Phi_\infty), \qquad \lim_{k\to\infty}\area(\vec \Phi_k) = \area(\vec \Phi_\infty).
\end{equation*} 
After the mentioned re-parametrisations, the $\vec \Phi_k$'s are weakly conformal which implies $\Delta_k \vec \Phi_k = 2H_{\vec{\Phi}_k}$ 
where in any local chart $x$, the intrinsic Laplacians $\Delta_k$ are given by 
\begin{equation*}
	\Delta_k\vec\Phi = \frac{1}{\sqrt{g_k}}\partial_{x^i}(\sqrt{g_k}(g_k)^{ij}\partial_{x^j} \vec\Phi), \qquad g_k = \det (g_k)_{ij}
\end{equation*}
for all $\vec\Phi\in W^{2,2}(\mathbb T^2,\mathbb R^3)$ where $(g_k)_{ij} = \partial_{x^i}\vec\Phi_k\cdot\partial_{x^j}\vec\Phi_k$ and $(g_k)^{ij}$ is the inverse of $(g_k)_{ij}$.
Therefore, by the weak convergence \eqref{eq:Helfrich-weak_convergence} it follows that for all $\delta>0$,
\begin{equation*}
	\lim_{k\to\infty}\int_{\mathbb T^2\setminus\bigcup_{i=1}^NB_\delta(a_i)}H_{\vec\Phi_k}\,\mathrm d\mu_{\vec\Phi_k} = \int_{\mathbb T^2\setminus\bigcup_{i=1}^NB_\delta(a_i)}H_{\vec\Phi_\infty}\,\mathrm d\mu_{\vec\Phi_\infty}.
\end{equation*}
Moreover, by \cite[Equation~(4.7)]{MR3176354},
\begin{equation*}
	\liminf_{\delta \to 0}\liminf_{k\to\infty} \int_{B_\delta(a_i)}1\,\mathrm d\mu_{\vec\Phi_k} = 0
\end{equation*}
for all $i\in\{1,\ldots,N\}$. Using the Cauchy--Schwarz inequality and the uniform bound on the Willmore energy \eqref{eq:Helfrich-8pi_bound}, it follows that after passing to a subsequence
\begin{equation*}
		\lim_{k\to\infty}\int_{\mathbb T^2}H_{\vec\Phi_k}\,\mathrm d\mu_{\vec\Phi_k} = \int_{\mathbb T^2}H_{\vec\Phi_\infty}\,\mathrm d\mu_{\vec\Phi_\infty}.
\end{equation*} 
Thus, by lower semi continuity of the Willmore functional under the convergence of \eqref{eq:Helfrich-weak_convergence},
\begin{equation*}
	\mathcal H_{c_0}(\vec\Phi_\infty) \leq \liminf_{k\to\infty}\mathcal H_{c_0}(\vec\Phi_k), \qquad \mathcal W(\vec\Phi_\infty) < 8\pi.
\end{equation*}
Therefore, $\vec\Phi_\infty$ is a minimiser and, by the Li--Yau inequality, $\vec\Phi_\infty \in W^{2,2}(\mathbb T^2,\mathbb R^3)$ is an embedding without branch points 
since branch points have multiplicity at least $2$ according to \cite[Theorem~3.1]{MR2928715}.
Moreover, by the regularity result \cite[Theorem~4.3]{MR4076069} (after~\cite{MR2430975}), $\vec \Phi_\infty \in C^\infty(\mathbb T^2,\mathbb R^3)$ which completes the proof of Corollary~\ref{cor:Helfrich-tori}. \qed


\begin{thebibliography}{DMSS20}
	
	\bibitem[BBM14]{BBM2014}
	Enrique Bendito, Mark~J. Bowick, and Agust\'{\i}n Medina.
	\newblock A natural parameterization of the roulettes of the conics generating
	the {D}elaunay surfaces.
	\newblock {\em J. Geom. Symmetry Phys.}, 33:27--45, 2014.
	
	\bibitem[BdC84]{MR731682}
	Jo\~{a}o~Lucas Barbosa and Manfredo do~Carmo.
	\newblock Stability of hypersurfaces with constant mean curvature.
	\newblock {\em Math. Z.}, 185(3):339--353, 1984.
	
	\bibitem[BF71]{eInts}
	Paul~F. Byrd and Morris~D. Friedman.
	\newblock {\em Handbook of elliptic integrals for engineers and scientists}.
	\newblock Die Grundlehren der mathematischen Wissenschaften, Band 67.
	Springer-Verlag, New York-Heidelberg, 1971.
	\newblock Second edition, revised.
	
	\bibitem[Bla09]{MR2591055}
	Simon Blatt.
	\newblock A singular example for the {W}illmore flow.
	\newblock {\em Analysis (Munich)}, 29(4):407--430, 2009.
	
	\bibitem[BLS20]{MR4098040}
	Katharina Brazda, Luca Lussardi, and Ulisse Stefanelli.
	\newblock Existence of varifold minimizers for the multiphase
	{C}anham-{H}elfrich functional.
	\newblock {\em Calc. Var. Partial Differential Equations}, 59(3):Paper No. 93,
	26, 2020.
	
	\bibitem[CV13]{MR3116014}
	Rustum Choksi and Marco Veneroni.
	\newblock Global minimizers for the doubly-constrained {H}elfrich energy: the
	axisymmetric case.
	\newblock {\em Calc. Var. Partial Differential Equations}, 48(3-4):337--366,
	2013.
	
	\bibitem[CVG07]{castro2007inverted}
	Pavel Castro-Villarreal and Jemal Guven.
	\newblock Inverted catenoid as a fluid membrane with two points pulled
	together.
	\newblock {\em Physical Review E}, 76(1):011922, 2007.
	
	\bibitem[Del41]{Delaunay1841}
	Charles~E. Delaunay.
	\newblock Sur la surface de r{\'e}volution dont la courbure moyenne est
	constante.
	\newblock {\em J. Math. Pures Appl.}, 6:309--314, 1841.
	
	\bibitem[DMSS20]{dall2020willmore}
	Anna Dall'Acqua, Marius M{\"u}ller, Reiner Sch{\"a}tzle, and Adrian Spener.
	\newblock The willmore flow of tori of revolution, 2020.
	
	\bibitem[Eic20]{MR4129521}
	Sascha Eichmann.
	\newblock Lower semicontinuity for the {H}elfrich problem.
	\newblock {\em Ann. Global Anal. Geom.}, 58(2):147--175, 2020.
	
	\bibitem[Hel73]{helfrich1973elastic}
	Wolfgang Helfrich.
	\newblock Elastic properties of lipid bilayers: theory and possible
	experiments.
	\newblock {\em Zeitschrift f{\"u}r Naturforschung C}, 28(11-12):693--703, 1973.
	
	\bibitem[H{\'e}l02]{MR1913803}
	Fr{\'e}d{\'e}ric H{\'e}lein.
	\newblock {\em Harmonic maps, conservation laws and moving frames}, volume 150
	of {\em Cambridge Tracts in Mathematics}.
	\newblock Cambridge University Press, Cambridge, second edition, 2002.
	\newblock Translated from the 1996 French original, With a foreword by James
	Eells.
	
	\bibitem[HMO07]{MR2381785}
	Mariana Hadzhilazova, Iva\"{\i}lo~M. Mladenov, and John Oprea.
	\newblock Unduloids and their geometry.
	\newblock {\em Arch. Math. (Brno)}, 43(5):417--429, 2007.
	
	\bibitem[HN21]{MR4270047}
	Lynn Heller and Cheikh~Birahim Ndiaye.
	\newblock First explicit constrained {W}illmore minimizers of non-rectangular
	conformal class.
	\newblock {\em Adv. Math.}, 386:Paper No. 107804, 47, 2021.
	
	\bibitem[IT92]{MR1215481}
	Y.~Imayoshi and M.~Taniguchi.
	\newblock {\em An introduction to {T}eichm\"{u}ller spaces}.
	\newblock Springer-Verlag, Tokyo, 1992.
	\newblock Translated and revised from the Japanese by the authors.
	
	\bibitem[Kap91]{MR1100207}
	Nicolaos Kapouleas.
	\newblock Compact constant mean curvature surfaces in {E}uclidean three-space.
	\newblock {\em J. Differential Geom.}, 33(3):683--715, 1991.
	
	\bibitem[Kap95]{MR1317648}
	Nikolaos Kapouleas.
	\newblock Constant mean curvature surfaces constructed by fusing {W}ente tori.
	\newblock {\em Invent. Math.}, 119(3):443--518, 1995.
	
	\bibitem[KKS89]{MR1010168}
	Nicholas~J. Korevaar, Rob Kusner, and Bruce Solomon.
	\newblock The structure of complete embedded surfaces with constant mean
	curvature.
	\newblock {\em J. Differential Geom.}, 30(2):465--503, 1989.
	
	\bibitem[KL12]{MR2928715}
	Ernst Kuwert and Yuxiang Li.
	\newblock {$W^{2,2}$}-conformal immersions of a closed {R}iemann surface into
	{$\mathbb R^n$}.
	\newblock {\em Comm. Anal. Geom.}, 20(2):313--340, 2012.
	
	\bibitem[KL18]{MR3842922}
	Ernst Kuwert and Yuxiang Li.
	\newblock Asymptotics of {W}illmore minimizers with prescribed small
	isoperimetric ratio.
	\newblock {\em SIAM J. Math. Anal.}, 50(4):4407--4425, 2018.
	
	\bibitem[KMR14]{MR3176354}
	Laura G.~A. Keller, Andrea Mondino, and Tristan Rivi\`ere.
	\newblock Embedded surfaces of arbitrary genus minimizing the {W}illmore energy
	under isoperimetric constraint.
	\newblock {\em Arch. Ration. Mech. Anal.}, 212(2):645--682, 2014.
	
	\bibitem[KP86]{MR868618}
	Wolfgang K\"{u}hnel and Ulrich Pinkall.
	\newblock On total mean curvatures.
	\newblock {\em Quart. J. Math. Oxford Ser. (2)}, 37(148):437--447, 1986.
	
	\bibitem[KS04]{MR2119722}
	Ernst Kuwert and Reiner Sch\"{a}tzle.
	\newblock Removability of point singularities of {W}illmore surfaces.
	\newblock {\em Ann. of Math. (2)}, 160(1):315--357, 2004.
	
	\bibitem[KS13]{MR3024303}
	Ernst Kuwert and Reiner Sch\"{a}tzle.
	\newblock Minimizers of the {W}illmore functional under fixed conformal class.
	\newblock {\em J. Differential Geom.}, 93(3):471--530, 2013.
	
	\bibitem[Kus89]{MR996204}
	Rob Kusner.
	\newblock Comparison surfaces for the {W}illmore problem.
	\newblock {\em Pacific J. Math.}, 138(2):317--345, 1989.
	
	\bibitem[Law70]{MR0270280}
	H.~Blaine Lawson, Jr.
	\newblock Complete minimal surfaces in {$S^{3}$}.
	\newblock {\em Ann. of Math. (2)}, 92:335--374, 1970.
	
	\bibitem[LY82]{MR674407}
	Peter Li and Shing~Tung Yau.
	\newblock A new conformal invariant and its applications to the {W}illmore
	conjecture and the first eigenvalue of compact surfaces.
	\newblock {\em Invent. Math.}, 69(2):269--291, 1982.
	
	\bibitem[MO03]{MR1977569}
	Iva\"{\i}lo Mladenov and John Oprea.
	\newblock Unduloids and their closed-geodesics.
	\newblock In {\em Geometry, integrability and quantization ({S}ts.
		{C}onstantine and {E}lena, 2002)}, pages 206--234. Coral Press Sci. Publ.,
	Sofia, 2003.
	
	\bibitem[MP01]{MR1807955}
	Rafe Mazzeo and Frank Pacard.
	\newblock Constant mean curvature surfaces with {D}elaunay ends.
	\newblock {\em Comm. Anal. Geom.}, 9(1):169--237, 2001.
	
	\bibitem[MR14]{MR3229052}
	Stefan M\"{u}ller and Matthias R\"{o}ger.
	\newblock Confined structures of least bending energy.
	\newblock {\em J. Differential Geom.}, 97(1):109--139, 2014.
	
	\bibitem[M{\v{S}}95]{MR1366547}
	Stefan M\"{u}ller and Vladim\'{\i}r {\v{S}}ver\'{a}k.
	\newblock On surfaces of finite total curvature.
	\newblock {\em J. Differential Geom.}, 42(2):229--258, 1995.
	
	\bibitem[MS20]{MR4076069}
	Andrea Mondino and Christian Scharrer.
	\newblock Existence and regularity of spheres minimising the
	{C}anham-{H}elfrich energy.
	\newblock {\em Arch. Ration. Mech. Anal.}, 236(3):1455--1485, 2020.
	
	\bibitem[MS21]{mondino2020strict}
	Andrea Mondino and Christian Scharrer.
	\newblock A strict inequality for the minimization of the {W}illmore functional
	under isoperimetric constraint.
	\newblock {\em Advances in Calculus of Variations}, 2021.
	
	\bibitem[NS14]{MR3247390}
	Cheikh~Birahim Ndiaye and Reiner~Michael Sch\"{a}tzle.
	\newblock Explicit conformally constrained {W}illmore minimizers in arbitrary
	codimension.
	\newblock {\em Calc. Var. Partial Differential Equations}, 51(1-2):291--314,
	2014.
	
	\bibitem[NS15]{MR3403429}
	Cheikh~Birahim Ndiaye and Reiner~Michael Sch\"{a}tzle.
	\newblock New examples of conformally constrained {W}illmore minimizers of
	explicit type.
	\newblock {\em Adv. Calc. Var.}, 8(4):291--319, 2015.
	
	\bibitem[Riv08]{MR2430975}
	Tristan Rivi\`ere.
	\newblock Analysis aspects of {W}illmore surfaces.
	\newblock {\em Invent. Math.}, 174(1):1--45, 2008.
	
	\bibitem[Riv13]{MR3008339}
	Tristan Rivi\`ere.
	\newblock Lipschitz conformal immersions from degenerating {R}iemann surfaces
	with {$L^2$}-bounded second fundamental forms.
	\newblock {\em Adv. Calc. Var.}, 6(1):1--31, 2013.
	
	\bibitem[Riv14]{MR3276154}
	Tristan Rivi\`ere.
	\newblock Variational principles for immersed surfaces with {$L^2$}-bounded
	second fundamental form.
	\newblock {\em J. Reine Angew. Math.}, 695:41--98, 2014.
	
	\bibitem[Riv15]{MR3383803}
	Tristan Rivi\`ere.
	\newblock Critical weak immersed surfaces within sub-manifolds of the
	{T}eichm\"uller space.
	\newblock {\em Adv. Math.}, 283:232--274, 2015.
	
	\bibitem[Sch12]{Schygulla}
	Johannes Schygulla.
	\newblock Willmore minimizers with prescribed isoperimetric ratio.
	\newblock {\em Arch. Ration. Mech. Anal.}, 203(3):901--941, 2012.
	
	\bibitem[Sch21]{scharrer:PhD}
	Christian Scharrer.
	\newblock On the minimisation of bending energies related to the willmore
	functional under constraints on area and volume, 2021.
	\newblock PhD thesis, \emph{Institutional Repository of the University of
		Warwick}, 2021.
	
	\bibitem[Sim93]{MR1243525}
	Leon Simon.
	\newblock Existence of surfaces minimizing the {W}illmore functional.
	\newblock {\em Comm. Anal. Geom.}, 1(2):281--326, 1993.
	
	\bibitem[Wil65]{MR0202066}
	Thomas~J. Willmore.
	\newblock Note on embedded surfaces.
	\newblock {\em An. \c{S}ti. Univ. ``Al. I. Cuza'' Ia\c{s}i Sec\c{t}. I a Mat.
		(N.S.)}, 11B:493--496, 1965.
	
	\bibitem[Woj17]{MR3744688}
	Stephan Wojtowytsch.
	\newblock Helfrich's energy and constrained minimisation.
	\newblock {\em Commun. Math. Sci.}, 15(8):2373--2386, 2017.
	
	\bibitem[YC22]{yu2020uniqueness}
	Thomas Yu and Jingmin Chen.
	\newblock Uniqueness of {C}lifford torus with prescribed isoperimetric ratio.
	\newblock {\em Proc. Amer. Math. Soc.}, 150(4):1749--1765, 2022.
	
\end{thebibliography}
\end{document}